\documentclass[letterpaper, paper,11pt]{AAS}		

\usepackage{bm}
\usepackage{amsmath}
\usepackage{amsfonts}
\usepackage[colorlinks=true, pdfstartview=FitV, linkcolor=black, citecolor= black, urlcolor= black]{hyperref}
\usepackage{overcite}
\usepackage[utf8]{inputenc}
\usepackage{textcomp}
\usepackage{pgfplots}
\usepackage{subcaption}
\usepackage{siunitx}
\usepackage{footnpag}	

\PaperNumber{24-440}

\begin{document}

\title{ADAPTIVE CONTROLLER FOR SIMULTANEOUS SPACECRAFT ATTITUDE TRACKING AND REACTION WHEEL FAULT DETECTION}

\author{Camilo Riano-Rios\thanks{Assistant Professor, Aerospace, Physics and Space Sciences Department, Florida Institute of Technology, 150 W. University Blvd., Melbourne, FL, 32940.},  
George Nehma\thanks{Ph.D. Student, Aerospace, Physics and Space Sciences Department, Florida Institute of Technology, 150 W. University Blvd., Melbourne, FL, 32940.},
\ and Madhur Tiwari\thanks{Assistant Professor, Aerospace, Physics and Space Sciences Department, Florida Institute of Technology, 150 W. University Blvd., Melbourne, FL, 32940.}
}

\maketitle{}

\begin{abstract}
The attitude control of a spacecraft is integral to achieving mission success. However, failures in actuators such as reaction wheels are detrimental and can often lead to an early end of mission. We propose a Lyapunov-based adaptive controller that can estimate and compensate for reaction wheels degradation simultaneously. The controller incorporates an adaptive update control law with a gradient-based term and an integral concurrent learning term that collects input-output data for online estimation of uncertain parameters. The proposed controller guarantees attitude tracking and its performance is tested through numerical simulations. 
\end{abstract}

\section{Introduction}

The spacecraft’s ability to accurately modify its attitude is one of the most critical features that must be guaranteed in any space mission. Reaction Wheels (RW) are the most commonly used actuators for this purpose due to their relatively small size and torque accuracy. However, their moving parts and the harsh conditions in space make these actuators prone to failure. The  Attitude Control System (ACS) of a spacecraft is a highly critical system, hence the failure of this system often means the end of the mission. Redundant Reaction Wheel Arrays (RWA) with more than three RWs are often used to maintain control authority for three-axis attitude control, and for momentum management using the RWA null space \cite{Evers2017}. 

Designing controllers for attitude control requires that they be resilient to disturbances and uncertainties. Current methods of control to address these issues include, sliding mode control \cite{Li2022-nt}, adaptive control \cite{Xiao2022-su}, observers \cite{Liang2019-ju} and neural networks \cite{rs13122396}. Often due to the lack of accuracy and difficulty in proving stability, neural networks are rejected as viable solutions, whilst observers  are used at the cost of more computational burden. Hence, of particular importance is the development of sliding mode control, which is a powerful yet simple controller that unfortunately can often suffer from chattering. Reference~\citenum{M_Sadigh2023-pl} develops an adaptive, fault tolerant sliding mode that does not generate chattering through an adaptive control protocol. Adaptive controllers are also a prominent solution because of their innate ability to handle uncertainties in the nonlinear dynamics. \cite{Xie2023-qc,M_Sadigh2023-pl,Wang2020-mi}

Fault-tolerant controllers have been developed to compensate for the reduction in performance of the ACS actuators. In Reference~\citenum{Fazlyab2023-cc}, a backstepping sliding mode controller was designed to maneuver the spacecraft’s attitude despite faulty actuators, which were a combination of RWs and thrusters. Reference~\citenum{Mei2023-xn} develops a sliding mode combined with disturbance observer fault tolerant controller that can quickly and efficiently handle actuator saturation and misalignment. Sensor failures were also detected and isolated using Unscented Kalman filters (UKF). Authors in Reference~\citenum{Chen2021-fe} designed a RW failure detection strategy that employs a pair of long short-term memory (LSTM) neural networks, one to approximate the RW dynamics and the other to train an adaptive threshold detector based on the difference between RW measurements and the first network.

A fault diagnosis strategy was designed in Reference~\citenum{Alidadi2023-cs} using a particle filter to predict the non-measurable bearing temperature lubricant in a RW, which has been identified as one of the main indicators for potential failures in this type of actuator. A similar approach was proposed in Reference~\citenum{Park2023-nr} to predict the remaining useful life of RWs by using an adaptive extended Kalman filter. In Reference~\citenum{Xiao2013-qz} a Lyapunov-based terminal sliding-mode observer was proposed to reconstruct RW faults and disturbances, and a compensation control law was designed to account for them. 

Adaptive control strategies are often used to compensate for parametric uncertainties in the system model by proposing gradient-based adaptation laws that compute estimates of the uncertain parameters online \cite{khalil2002nonlinear}. These adaptation laws are often sufficient to prove convergence of the error states. However, convergence of the model parameters to their real values can only be guaranteed under the assumption of Persistent Excitation (PE), i.e., the system is persistently excited over an infinite time integral, which is difficult to guarantee and verify in practice \cite{ioannou2012robust, sastry2011adaptive}. An adaptation strategy, called Concurrent Learning (CL), was proposed to relax the PE condition by introducing in the adaptation law a term that collects input and output data under the assumption of measurable higher order states \cite{Chowdhary2013}. The CL adaptation strategy guarantees exponential convergence of the system parameter estimates after a verifiable Finite Excitation (FE) condition is satisfied, which can be evaluated online, and concurrent to the control execution, by computing the eigenvalue of a matrix built with collected input-output data. Recent efforts have relaxed the assumption of measurable higher order states by using their integrals in the so-called Integral Concurrent Learning (ICL) \cite{kamalapurkar2019}, which has been used in a wide range of control applications \cite{Bell2020,RIANORIOS2020189,Sun2021}. 

Compensation for actuator fault detection and detection or monitoring of actuator degradation are often considered as separate blocks in the spacecraft ACS, in the control and estimation blocks, respectively. This paper introduces a Lyapunov-based adaptive controller that performs both compensation and fault detection simultaneously, and achieves global exponential attitude tracking in the presence of RW failure or degradation. An actuator health matrix is included in the satellite attitude dynamics and considered uncertain. An ICL-based adaptive update law is designed to achieve convergence of the RW health parameters once the FE condition is satisfied. The controller validated in simulation with different number of RWs and degradation levels. 

The paper is organized as follows: A first section presents the satellite attitude dynamic model, the next two sections describe the controller design and corresponding stability analysis. Finally, sections with preformed numerical simulations and conclusions are presented.

\section{Spacecraft Attitude Dynamics}

\begin{figure}[htb]
	\centering\includegraphics[width=3.5in]{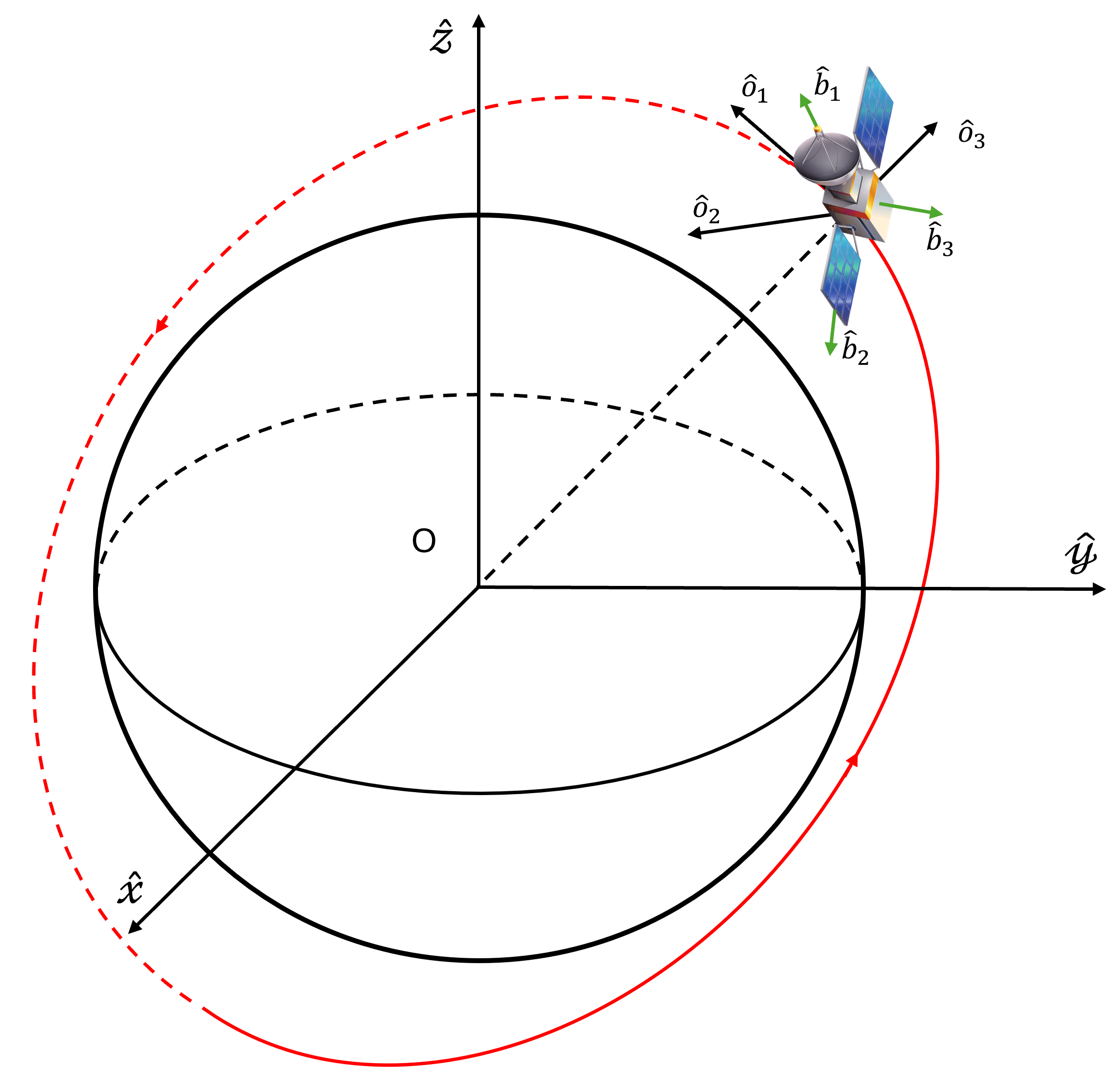}
	\caption{Coordinate Systems}
	\label{fig:coord_systems}
\end{figure}

Three reference frames and their corresponding coordinate systems will be considered throughout this paper, and are depicted in Figure \ref{fig:coord_systems}. The Body reference frame $\mathcal{B}$, attached to the spacecraft body is defined with origin at its center of mass, and the unit vectors $\boldsymbol{\hat{b}_1, \hat{b}_2}$ and $\boldsymbol{\hat{b}_3}\in\mathbb{R}^3$ along its longitudinal, lateral and vertical axes, respectively. The coordinate system associated with the orbital reference frame $\mathcal{O}$ has origin at the spacecraft center of mass, and is defined by the unit vectors $\boldsymbol{\hat{o}_1, \hat{o}_2}$ and $\boldsymbol{\hat{o}_3}\in\mathbb{R}^3$, with $\boldsymbol{\hat{o}_3}$ aligned with the zenith direction, $\boldsymbol{\hat{o}_2}$ the orbital angular momentum vector, and $\boldsymbol{\hat{o}_1}$ aligned with the $(\boldsymbol{\hat{o}_2 \times \hat{o}_3})$ direction. The Earth-Centered-Inertial (ECI) coordinate system, centered at the Earth's center of mass and with basis $\{\boldsymbol{\hat{x}, \hat{y},\hat{z}}\}$, is associated with the Inertial reference frame $\mathcal{N}$.

\subsection{Equation of Motion}

The attitude motion of a spacecraft with $N\in\mathbb{Z}_{>0}$ reaction wheels is governed by the following equations of motion:

\begin{align}
    \label{eq:Eulers_law}
	J\boldsymbol{\dot{\omega}} &= -\boldsymbol{\omega} \times \left(J\boldsymbol{\omega} + J_{RW}G\boldsymbol{\Omega} \right) + G\Phi\boldsymbol{u}  \\ 
        \boldsymbol{\dot{\sigma}} &= \frac{1}{4}\left[\left( 
1-\boldsymbol{\sigma^T\sigma}\right)I_3 + 2\sigma^{\times}+2\boldsymbol{\sigma\sigma^T}\right]\boldsymbol{\omega},\label{eq:MRP_kinem}
\end{align}
where $\boldsymbol{\omega}\in\mathbb{R}^3$ is the spacecraft angular velocity expressed in the body coordinate system, $\boldsymbol{\sigma\in\mathbb{R}^3}$ is the vector of Modified Rodrigues Parameters (MRP) that represent the orientation of the spacecraft with respect to the inertial frame\cite{Schaub2018}, $\boldsymbol{\Omega}=[\Omega_1, \Omega_2, \cdots, \Omega_N]^T~\in\mathbb{R}^N$ is a vector containing the $N$ RW angular velocities, $\boldsymbol{u}=-J_{RW}\boldsymbol{\dot{\Omega}}=[u_1, u_2, \cdots, u_N]^T ~ \in \mathbb{R}^N$ is the control input that represents the torque applied by each RW, $J\in\mathbb{R}^{3\times 3}$ is the total inertia matrix, $J_{RW}\in\mathbb{R}_{>0}$ is the inertia of the flywheels about their spin axis, $\Phi=diag\{\phi_1,~\phi_2,~\cdots,~\phi_N\}\in\mathbb{R}^{N
\times N}$ is the uncertain RW health matrix, $G=\{\boldsymbol{\hat{s}_1, \hat{s}_2, \cdots, \hat{s}_N}\}\in\mathbb{R}^{3\times N}$ is the RWA configuration matrix, and $\boldsymbol{\hat{s}_i}\in\mathbb{R}^3$ is the direction of the $i^{th}$ RW's spin axis expressed in the body coordinate system. The matrix $I_m\in\mathbb{R}^{m\times m}$ represents an identity matrix of dimension $m\times m$, and the skew-symmetric matrix $a^{\times}\in\mathbb{R}^{3\times 3}$ for a vector $\boldsymbol{a}=[a_1, a_2, a_3]^T\in\mathbb{R}^3$ is defined as 

\begin{equation}
    \label{eq:skew_sym_mat}
    a^{\times}=\begin{bmatrix}
    0 & -a_3 & a_2 \\ 
    a_3 & 0 & -a_1 \\
    -a_2 & a_1 & 0
    \end{bmatrix}.
\end{equation}

\section{Control Design}
    
\subsection{Control Objective}

The objective is to design an adaptive controller that achieves attitude tracking despite failures in one or more redundant reaction wheels. Besides compensating for these failures, the controller must also provide an online estimate of the health matrix $\Phi$ for actuators performance assessment.   

\paragraph{Assumption 1.} Changes or degradation in performance of RWs are slow enough, as compared to attitude maneuver times, for the $\Phi$ matrix to be considered constant.  

The spacecraft attitude dynamics in Eq. \ref{eq:Eulers_law} can be rewritten as

\begin{equation}
    \label{eq:EoMs_Ytheta}
    \boldsymbol{\dot{\omega}} = J^{-1} \left(-\boldsymbol{\omega} \times \left(J\boldsymbol{\omega}+J_{RW}G\boldsymbol{\Omega}\right) + Y\boldsymbol{\theta}\right),
\end{equation}
where the term $Y\boldsymbol{\theta}=G\Phi\boldsymbol{u} ~ \in\mathbb{R}^3$ is a linear parameterization with respect to the uncertain diagonal entries of the health matrix $\Phi$. $Y\in\mathbb{R}^{3\times N}$ is a measurable regression matrix defined as 

\begin{equation}
    \label{eq:regress_matrix}
    Y = G\left(diag\{u_1, u_2, \cdots, u_N\}\right),
\end{equation}
and the vector of uncertain parameters $\boldsymbol{\theta}\in\mathbb{R}^{N}$ as

\begin{equation}
    \label{eq:theta}
    \boldsymbol{\theta}=\left[\phi_1,~\phi_2,~\cdots,~\phi_N\right]^T
\end{equation}

\paragraph{Assumption 2.} The spacecraft is equipped with an attitude determination system capable of providing the controller with angular velocity $\boldsymbol{\omega}$ and attitude $\boldsymbol{\sigma}$ measurements.

Let us introduce the error MRP $\boldsymbol{\sigma_e}\in\mathbb{R}^3$, representing the attitude mismatch between the body and desired frames, which can be computed using the spacecraft attitude $\boldsymbol{\sigma}$ and desired, bounded attitude trajectories $\boldsymbol{\sigma_d,\omega_d}\in\mathbb{R}^3$, and obeys the following kinematic equation\cite{Schaub2018}

\begin{equation}
    \label{eq:error_mrp_kinem}
    \boldsymbol{\dot{\sigma}_e}=\frac{1}{4}B\boldsymbol{\Tilde{\omega}},
\end{equation}
where the matrix $B\in\mathbb{R}^{3\times 3}$ is defined as 

\begin{equation}
    \label{eq:B_matrix}
    B = \left(1-\boldsymbol{\sigma_e^T \sigma_e}\right)I_3 + 2\sigma_e^{\times} + 2\boldsymbol{\sigma_e\sigma_e^T}, 
\end{equation}
 $\boldsymbol{\Tilde{\omega}}=\boldsymbol{\omega}-\Tilde{R}\boldsymbol{\omega}_d~\in\mathbb{R}^3$ is the relative angular velocity and $\Tilde{R}\in\mathbb{R}^{3\times 3}$, defined as 

\begin{equation}
    \label{eq:R_tilde}
    \Tilde{R}=I_3+\frac{8\left(\boldsymbol{\sigma_e^{\times}}\right)-4\left(1-\boldsymbol{\sigma_e^T\sigma_e}\right)\boldsymbol{\sigma_e^{\times}}}{\left(1+\boldsymbol{\sigma_e^T\sigma_e}\right)^2},
\end{equation}
represents the rotation matrix between the body and desired frames. 

The attitude control objective can be established as

\begin{equation}
    \label{eq:ctrl_obj1}
    \Tilde{R}\rightarrow I_3, ~~~ \mathrm{as}~~t\rightarrow\infty,  
\end{equation}
which will be achieved if 

\begin{equation}
    \label{eq:ctrl_obj2}
    \left\|\boldsymbol{\sigma_e}\right\|\rightarrow 0,~~ \mathrm{and} ~~ \left\|\boldsymbol{\Tilde{\omega}}\right\|\rightarrow 0 ~\Rightarrow~\left\|\boldsymbol{\dot{\sigma}_e}\right\|\rightarrow 0.
\end{equation}

\subsection{Control Development}

Let us define a modified state vector $\boldsymbol{r}\in\mathbb{R}^3$ as 

\begin{equation}
    \label{eq:state_r}
    \boldsymbol{r}=\boldsymbol{\dot{\sigma_e}}+\alpha\boldsymbol{\sigma_e},
\end{equation}
where $\alpha\in\mathbb{R}^{3\times 3}$ is a symmetric, positive definite control gain matrix. Taking the time derivative of Eq. \ref{eq:state_r} yields

\begin{equation}
        \label{eq:open_loop_error_sys}
        \boldsymbol{\dot{r}}=\frac{1}{4} \dot{B}\boldsymbol{\Tilde{\omega}} + \frac{1}{4}B\left(\boldsymbol{\dot{\omega}}-\Tilde{R}\boldsymbol{\dot{\omega}_d}-\dot{\Tilde{R}}\boldsymbol{\omega_d}\right)+\alpha\boldsymbol{\dot{\sigma_e}},
\end{equation}    
with $\dot{B}=\left[-2\boldsymbol{\sigma_e^T\dot{\sigma}_e}I_3+2\dot{\sigma}_e^{\times}+4\left(\boldsymbol{\dot{\sigma}_e\sigma_e^T}\right)\right]$.

Using the dynamics in Eq. \ref{eq:EoMs_Ytheta}, and the fact that $\dot{\Tilde{R}}=-\Tilde{\omega}^{\times}\Tilde{R}$, we obtain 

\begin{equation}
    \label{eq:r_dot_interm}
    \boldsymbol{\dot{r}}=\frac{1}{4}\dot{B}\boldsymbol{\Tilde{\omega}}+\frac{1}{4}B\left(J^{-1}\left(-\boldsymbol{\omega}\times \left(J\boldsymbol{\omega}+J_{RW}G\boldsymbol{\Omega}\right)+Y\boldsymbol{\Tilde{\theta}}+Y\boldsymbol{\hat{\theta}}\right)-\Tilde{R}\boldsymbol{\dot{\omega}_d}+\Tilde{\omega}^{\times}\Tilde{R}\boldsymbol{\omega_d}\right)+\alpha\boldsymbol{\dot{\sigma}_e},
\end{equation}
where the estimation error is defined as $\boldsymbol{\Tilde{\theta}}=\boldsymbol{\theta-\hat{\theta}}~\in\mathbb{R}^{N}$, and therefore the fact that $Y\boldsymbol{\theta}=Y\boldsymbol{\tilde{\theta}}+Y\boldsymbol{\hat{\theta}}$ was used. 

Since the control input $\boldsymbol{u}$, i,e., RW torques, is embedded in the $Y\boldsymbol{\theta}$ term, let the auxiliary control signal $\boldsymbol{u_d}\in\mathbb{R}^3$ be designed as

\begin{equation}
    \label{eq:ctrl_law}
    \boldsymbol{u_d}=\boldsymbol{\omega}\times \left(J\boldsymbol{\omega}+J_{RW}G\boldsymbol{\Omega}\right) + J\Tilde{R}\boldsymbol{\dot{\omega}_d}-J\Tilde{\omega}^{\times}\Tilde{R}\boldsymbol{\omega_d}+4JB^{-1}\left[-\frac{1}{4}\dot{B}\boldsymbol{\tilde{\omega}}-\alpha\boldsymbol{\dot{\sigma}_e}-K\boldsymbol{r}-\beta\boldsymbol{\sigma_e}\right],
\end{equation}
where $\beta\in\mathbb{R}_{>0}$ is a constant control gain. 

Based on Eqs. \ref{eq:r_dot_interm} and \ref{eq:ctrl_law}, the adaptation law is designed as 

\begin{equation}
    \label{eq:adapt_law}
    \boldsymbol{\dot{\hat{\theta}}}=\mathrm{proj}\left\{\frac{1}{4}\Gamma Y^T\left(J^{-1}\right)^TB^T\boldsymbol{r}+\Gamma K_1\sum_{i=1}^{N_s}\mathcal{Y}_i^T\left(J\boldsymbol{\omega}(t)-J\boldsymbol{\omega}(t-\Delta t)+\boldsymbol{\mathcal{U}_i}-\mathcal{Y}_i\boldsymbol{\hat{\theta}}\right)\right\},
\end{equation}
where $\Gamma,~K_1,~\in\mathbb{R}^{N\times N}$ are constant, positive definite adaptation gain matrices, $\boldsymbol{\hat{\theta}}\in\mathbb{R}^{N}$ is the estimate of $\boldsymbol{\theta}$, $\Delta t\in\mathbb{R}_{>0}$ is the time between samples, $\mathrm{proj}\{\cdot\}$ represents a projection algorithm to keep $\boldsymbol{\hat{\theta}}$ within known, user-defined bounds \cite{dixon2003nonlinear}, $N_s\in\mathbb{Z}_{>0}$ is the number of input-output data pairs used for online parameter estimation, and the terms $\boldsymbol{\mathcal{U}_i}\in\mathbb{R}^3$ and $\mathcal{Y}_i\in \mathbb{R}^{3\times N}$ are defined as 

\begin{align}
    \label{eq:ICL_terms}
    \boldsymbol{\mathcal{U}_i}(\Delta t, t_i)&=\int_{t_i-\Delta t}^{t_i}\left(\boldsymbol{\omega}(\tau)\times \left(J\boldsymbol{\omega}(\tau)+J_{RW}G\boldsymbol{\Omega}(\tau)\right)\right)d\tau, \\
    \mathcal{Y}_i(\Delta t, t_i)&=\int_{t_i-\Delta t}^{t_i}Y(\tau)d\tau.
\end{align}

Based on the definition of $Y\boldsymbol{\theta}$, the term $Y\boldsymbol{\hat{\theta}}$ can be written as $Y\boldsymbol{\hat{\theta}}=G\hat{\Phi} \boldsymbol{u}$, where $\hat{\Phi}$ is the estimate of $\Phi$. Note that the actual control input is $\boldsymbol{u}$. Then, we can set 

\begin{equation}
    \label{eq:Y_thetahat}
    Y\boldsymbol{\hat{\theta}}=\boldsymbol{u_d},
\end{equation}
and recover $\boldsymbol{u}$ as 

\begin{equation}
    \label{eq:u_real} \boldsymbol{u}=\left(G\hat{\Phi}\right)^{\dag}\boldsymbol{u_d}, 
\end{equation}
where $\hat{\Phi}$ can be obtained by numerically integrating Eq. \ref{eq:adapt_law}, and $(\cdot)^{\dag}$ is the pseudo-inverse of $(\cdot)$.

\paragraph{Assumption 3.} The system is sufficiently excited over a finite duration of time. Then, there exists a finite time $T\in\mathbb{R}_{>0}$ such that 

\begin{equation}
    \label{eq:FE_condition}
    \lambda_{min}\left\{\sum_{i=1}^{N_s}\mathcal{Y}_i^T\mathcal{Y}_i\right\}\geq\Bar{\lambda},
\end{equation}
where $\lambda_{min}\{\cdot\}$ denotes the minimum eigenvalue of the finite excitation condition matrix $\{\cdot\}$, and $\Bar{\lambda}\in\mathbb{R}_{>0}$ is a user-defined threshold.

\section{Stability Analysis}

For the subsequent stability analysis, let us define the composite state vector $\boldsymbol{\eta}=[\boldsymbol{r^T,~\sigma_e^T,~\Tilde{\theta}^T}]^T~\in\mathbb{R}^{6+N}$, and divide it into two theorems: the first theorem considers the stability of the closed-loop system before the finite excitation condition in Assumption 3 is satisfied, i.e., $t<T$; and the second theorem presents the stability of the closed-loop system after t=T.  

\paragraph{Theorem 1.} Given the spacecraft attitude dynamics in Eqs. \ref{eq:Eulers_law} and \ref{eq:MRP_kinem}, along with the adaptive update law in Eq. \ref{eq:adapt_law}, the controller proposed in Eq. \ref{eq:ctrl_law} ensures that the estimation error $\boldsymbol{\Tilde{\theta}}$ remains bounded and asymptotic attitude tracking is achieved so that conditions in Eqs. \ref{eq:ctrl_obj1} and \ref{eq:ctrl_obj2} are satisfied. 

\paragraph{Proof.} Let $V:\mathbb{R}^{6+N}\rightarrow \mathbb{R}_{>0}$ be a Lyapunov candidate function defined as 

\begin{equation}
    \label{eq:Lyap_func}
    V(\boldsymbol{\eta})=\frac{1}{2}\boldsymbol{r^Tr}+\frac{\beta}{2}\boldsymbol{\sigma_e^T\sigma_e}+\frac{1}{2}\boldsymbol{\Tilde{\theta}}^T\Gamma^{-1}\boldsymbol{\Tilde{\theta}}
\end{equation}

Taking the time derivative of Eq. \ref{eq:Lyap_func}, substituting Eqs. \ref{eq:state_r}, \ref{eq:ctrl_law} and \ref{eq:Y_thetahat}, we get

\begin{equation}
    \label{eq:V_dot_interm}
    \dot{V}(\boldsymbol{\eta})=\boldsymbol{r}^T\left(\frac{1}{4}BJ^{-1}Y\boldsymbol{\Tilde{\theta}}-K\boldsymbol{r}\right)-\beta\boldsymbol{\sigma_e}^T\alpha\boldsymbol{\sigma_e}-\boldsymbol{\Tilde{\theta}}^T\Gamma^{-1}\boldsymbol{\dot{\hat{\theta}}}.
\end{equation}

The expression in Eq. \ref{eq:adapt_law} can be represented in the non-implementable -but useful for the stability analysis- form 

\begin{equation}
    \label{eq:adapt_law_non_impl}
    \boldsymbol{\dot{\hat{\theta}}}=\mathrm{proj}\left\{\frac{1}{4}\Gamma Y^T\left(J^{-1}\right)^TB^T\boldsymbol{r}+\Gamma K_1\sum_{i=1}^{N_s}\mathcal{Y}_i^T\mathcal{Y}_i\boldsymbol{\Tilde{\theta}}\right\},
\end{equation}
which can be plugged in Eq. \ref{eq:V_dot_interm} to obtain 

\begin{equation}
    \label{eq:V_dot}
    \dot{V}(\boldsymbol{\eta})=-\boldsymbol{r}^TK\boldsymbol{r}-\beta\boldsymbol{\sigma_e}^T\alpha\boldsymbol{\sigma_e}-\boldsymbol{\Tilde{\theta}}^TK_1\sum_{i=1}^{N_s}\mathcal{Y}_i^T\mathcal{Y}_i\boldsymbol{\Tilde{\theta}}.
\end{equation}

Under Assumption 3, and before $t=T$, Eq. \ref{eq:V_dot} can be upper bounded as 

\begin{equation}
\label{eq:V_dot_upper_bound1}
    \dot{V}(\boldsymbol{\eta})\leq -\lambda_{min}\left\{K\right\}\left\|\boldsymbol{r}\right\|^2-\beta\lambda_{min}\left\{\alpha\right\}\left\|\boldsymbol{\sigma_e}\right\|^2.
\end{equation}

From Eq. \ref{eq:V_dot_upper_bound1}, $\dot{V}\in\mathcal{L}_{\infty}$. Then, $\boldsymbol{\eta}\in\mathcal{L}_{\infty}$ and $\boldsymbol{r,\sigma_e}\in\mathcal{L}_2$. Since $\boldsymbol{\sigma_e,r}\in\mathcal{L}_{\infty}$, then $\boldsymbol{\dot{\sigma}_e}\in\mathcal{L}_{\infty}$, which along with Eqs. \ref{eq:error_mrp_kinem} and \ref{eq:B_matrix}, and the fact that $\boldsymbol{\omega_d, \dot{\omega}_d},\Tilde{R}\in\mathcal{L}_{\infty}$ by definition, leads to $\boldsymbol{\Tilde{\omega}}\in\mathcal{L}_{\infty}\Rightarrow\boldsymbol{u_d}\in\mathcal{L}_{\infty}$. Since $\boldsymbol{u_d}\in\mathcal{L}_{\infty}$ and $\boldsymbol{\hat{\theta}}\in\mathcal{L}_{\infty}$ due to the projection algorithm, then $Y\in\mathcal{L}_{\infty}$, which implies that $\boldsymbol{\Omega},\boldsymbol{u}\in\mathcal{L}_{\infty}$, and therefore $\boldsymbol{\dot{r}}\in\mathcal{L}_{\infty}$ from Eq. \ref{eq:r_dot_interm}. Since $\boldsymbol{\dot{\sigma}_e,\dot{r}}\in\mathcal{L}_{\infty}$ and $\boldsymbol{r,\sigma_e}\in\mathcal{L}_2$, then by Barbalat's lemma \cite{khalil2002nonlinear}, we can conclude that attitude tracking is achieved. 

The subsequent theorem, assumes that the finite excitation condition in Eq. \ref{eq:FE_condition} is satisfied, which will allow online estimation of the RW health matrix $\Phi$. The purpose of estimating $\Phi$ is to obtain information about the RWs performance. For Instance, in the case of total failure on the $i^{th}$ RW the spacecraft will not experience torques coming from that RW. This will force the estimation of $\Phi's$ $i^{th}$ diagonal entry to approach zero provided guaranteed convergence of $\boldsymbol{\Tilde{\theta}}$, so that no torques from this RW influence the attitude dynamics. 

For the subsequent stability analysis let us define $\underline{\kappa},~\Bar{\kappa}~\in\mathbb{R}_{>0}$ as known constants used to bound the Lyapunov function, and $\zeta=\mathrm{min}\left(\lambda_{min}\left\{K\right\},\beta\lambda_{min}\left\{\alpha\right\},\lambda_{min}\left\{K_1\sum_{i=1}^{N_s}\mathcal{Y}_i^T\mathcal{Y}_i\right\}\right)~\in\mathbb{R}_{>0}$.

\paragraph{Theorem 2.} For the spacecraft attitude dynamics described in Eqs. \ref{eq:Eulers_law} and \ref{eq:MRP_kinem}, the controller in Eq. \ref{eq:ctrl_law} and the adaptation law in Eq. \ref{eq:adapt_law} achieve global exponential attitude tracking and convergence of the parameter estimation in the sense  that 

\begin{equation}
    \label{eq:result_thm_2}
    \left\|\boldsymbol{\eta}(t)\right\|\leq\frac{\Bar{\kappa}}{\underline{\kappa}}\mathrm{exp}\left(\frac{\zeta T}{2\Bar{\kappa}}\right)\left\|\boldsymbol{\eta}(0)\right\|\mathrm{exp}\left(-\frac{\zeta}{2\Bar{\kappa}}~t\right)~~\forall ~t\geq 0
\end{equation}

\paragraph{Proof.} The bounds for the Lyapunov function in Eq. \ref{eq:Lyap_func} can be expressed as 

\begin{equation}
    \label{eq:bounds_V}
    \underline{\kappa}\left\|\boldsymbol{\eta}\right\|^2\leq V(\boldsymbol{\eta})\leq\Bar{\kappa}\left\|\boldsymbol{\eta}\right\|^2,
\end{equation}
Since the matrix $\sum_{i=1}^{N_s}\mathcal{Y}_i^T\mathcal{Y}_i$ is positive definite after $t=T$, i.e., the finite excitation condition is satisfied, then the upper bound for $\dot{V}$ can be rewritten as 

\begin{equation}
    \label{eq:V_dot_upper_bnd_thm2}
    \dot{V}(\boldsymbol{\eta})\leq -\lambda_{min}\left\{K\right\}\left\|\boldsymbol{r}\right\|^2-\beta\lambda_{min}\left\{\alpha\right\}\left\|\boldsymbol{\sigma_e}\right\|^2-\lambda_{min}\left\{K_1\sum_{i=1}^{N_s}\mathcal{Y}_i^T\mathcal{Y}_i\right\}\left\|\boldsymbol{\tilde{\theta}}\right\|^2,
\end{equation}
or in a more compact form

\begin{equation}
    \label{eq:V_dot_upper_bnd_thm2_compact}
    \dot{V}(\boldsymbol{\eta})\leq-\zeta\left\|\boldsymbol{\eta}\right\|^2.
\end{equation} 

Applying the comparison lemma\cite{khalil2002nonlinear} and using the bounds in Eq. \ref{eq:bounds_V} yields

\begin{equation}
    \label{eq:stab_result1}
    V\left(\boldsymbol{\eta}(t)\right)\leq V(\boldsymbol{\eta}(T))\mathrm{exp}\left(-\frac{\zeta}{\Bar{\kappa}}\left(t-T\right)\right)~~\forall ~t\geq T.
\end{equation}
Note that in this portion of the analysis, the dependence of time of the state $\boldsymbol{\eta}(t)$ is made explicit. Using the bounds in Eq. \ref{eq:bounds_V} to express Eq. \ref{eq:stab_result1} in terms of $\left\|\boldsymbol{\eta}(t)\right\|$ and $\left\|\boldsymbol{\eta}(T)\right\|$ we get

\begin{equation}
    \label{eq:stab_result2}
    \left\|\boldsymbol{\eta}(t)\right\|\leq\sqrt{\frac{\Bar{\kappa}}{\underline{\kappa}}}\left\|\boldsymbol{\eta}(T)\right\|\mathrm{exp}\left(-\frac{\zeta}{2\Bar{\kappa}}(t-T)\right)~~\forall ~t\geq T.
\end{equation}
Note also that from Eq. \ref{eq:V_dot_upper_bnd_thm2_compact} we know that $\left\|\boldsymbol{\eta}(T)\right\|\leq\left\|\boldsymbol{\eta}(0)\right\|$. Using this fact and the bounds in Eq. \ref{eq:bounds_V} yields

\begin{equation}
    \label{eq:stab_result3}
    \left\|\boldsymbol{\eta}(t)\right\|\leq\frac{\Bar{\kappa}}{\underline{\kappa}}\left\|\boldsymbol{\eta}(0)\right\|\mathrm{exp}\left(-\frac{\zeta}{2\Bar{\kappa}}(t-T)\right)~~\forall ~t\geq 0, 
\end{equation}
which can finally be rewritten as shown in Eq. \ref{eq:result_thm_2}. 

\section{Simulation}

In order to demonstrate the capability of the designed controller to estimate the uncertain health parameters of the RWs, we devise 4 separate scenarios. In all simulation scenarios, the same attitude guidance commands are used to excite the system. In each case, the initial guidance is to align the satellite with the ECI frame. Then the guidance alternates between the initial configuration and Nadir pointing, i.e., alignment with the orbital frame, three times before settling with Nadir pointing. The change of attitude reference occurs every 12 minutes with the final Nadir pointing being set after 2000 seconds of simulation. The first two scenarios, both with 4 RWs aim to demonstrate the effect of the ICL term in the adaptation law. This is done by activating the ICL term (once the FE condition is satisfied) in the first simulation, and then repeating the same simulation with the ICL term disabled. 

The next pair of simulations illustrates the capability of the controller to identify the degradation level of a certain pair of RWs. Six RWs are used in both simulations and both contain two failing RWs. In the first scenario, the RWs are degraded to 0\% performance, i.e., complete failure, whilst in the second, RW 1 is degraded to 0\% performance and RW 2 is degraded to 30\% performance. The gains used for all scenarios are listed in Table \ref{tab:gains} whilst the parameters for the configuration of the RWs are listed in Table \ref{tab:config}.

\begin{table}[h]
	\fontsize{10}{10}\selectfont
    \caption{Table of gains for the adaptive controller for all 4 Cases.}
   \label{tab:gains}
        \centering
   \renewcommand{\arraystretch}{1.5}
   \begin{tabular}{c | c | c | c | c} 
      \hline 
      Gain    & Case 1 & Case 2 & Case 3 & Case 4\\
      \hline 
      $K_1$ & $10~I_4$ & $0$ & $10~I_6$ & $10~I_6$ \\
      $K$       & $5\times 10^{-1}~I_3$ & $5\times 10^{-1}~I_3$ & $5\times 10^{-2}~I_3$ & $5\times 10^{-2}~I_3$ \\
      $\alpha$  & $3\times 10^{-2}~I_3$ & $3\times 10^{-2}~I_3$ & $3\times 10^{-2}~I_3$ & $3\times 10^{-2}~I_3$ \\
      $\beta$       & $5\times 10^{-3}$ & $5\times 10^{-3}$ & $5\times 10^{-3}$ & $5\times 10^{-3}$ \\
      $\Gamma$ & $100~I_4$ & $100~I_4$ & $300~I_6$ & $300~I_6$ \\
      $\Bar{\lambda}$ & $1\times 10^{-7}$ & N/A & $8\times 10^{-9}$ & $8\times 10^{-9}$ \\
      $\Phi$ & $diag\{1,1,0,1\}$ & $diag\{1,1,0,1\}$ & $diag\{0,0,1,1,1,1\}$ & $diag\{0,0.3,1,1,1,1\}$ \\
      \hline
   \end{tabular}
\end{table}

\begin{table}[h]
	\fontsize{10}{10}\selectfont
    \caption{Table of configuration parameters for the RWs for each simulation.}
   \label{tab:config}
        \centering 
   \renewcommand{\arraystretch}{1.5}
   \begin{tabular}{c | c | c} 
      \hline 
      Parameter & Value & Units\\
      \hline 
      Satellite Mass & 25 & $\si{kg}$ \\
      $J$  & $diag\{0.4333, 0.7042, 0.7042\}$ & $\si{kg.m^2}$\\
      $G$ (4 RWs) & $\begin{bmatrix}
          0.5774 & -0.5774 & 0.5774 & -0.5774 \\
          0.5774 & 0.5774 & -0.5774 & -0.5774 \\
          0.5774 & 0.5774 & 0.5774 & 0.5774 \\
      \end{bmatrix}$ & N/A \\
      $G$ (6 RWs) & $\begin{bmatrix}
          0.5 & 0.5 & 0.5 & 0.5 & 0.5 & 0.5 \\
          0 & 0.75 & 0.75 & 0 & -0.75 & -0.75 \\
          0.866 & 0.433 & -0.433 & -0.866 & -0.433 & 0.433 \\
      \end{bmatrix}$ & N/A \\
      $J_{RW}$ & $5.7296\times 10^5$ & $\si{kg.m^2}$ \\
      Max RW Torque & $20\times 10^{-3}$ & $\si{N.m}$\\
      Max $\Omega$ (for all RWs) & $1.0472\times 10^3$ & $\si{rad.s^{-1}}$ \\
      \hline
   \end{tabular}
\end{table}

\subsection{Cases 1 \& 2: ICL term influence}

The ICL term is responsible for the input-output data collection and helps achieving convergence of the uncertain parameters' estimation error. Hence, the following two simulations illustrate the importance of this term in the controller's ability to estimate the health of the system's actuators; in this case, the health of all the RWs. 

\begin{figure}[h]
    \begin{subfigure}{0.45\textwidth} 
        \centering
        \includegraphics[width=\linewidth]{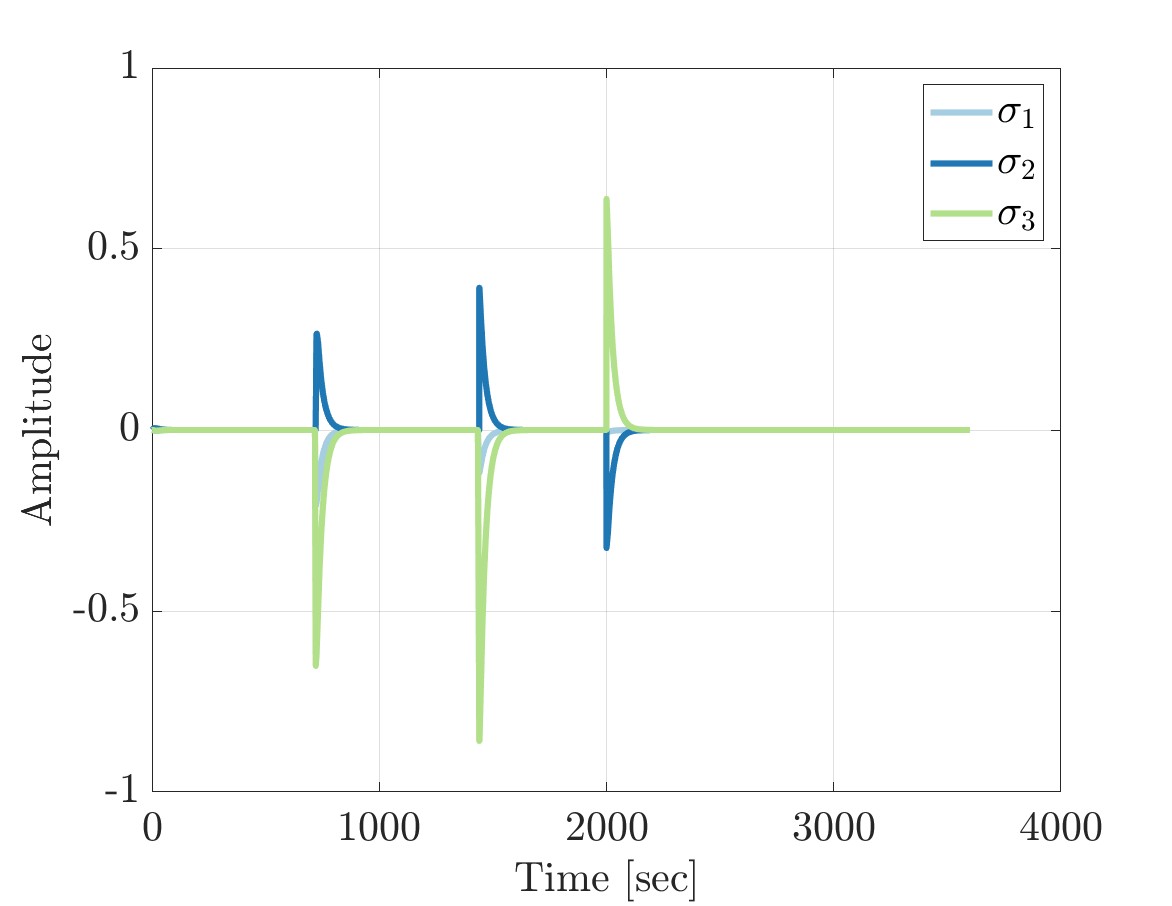}  
        \caption{Error MRP}
        \label{fig:RW4_ICL_MRP}
    \end{subfigure}
    \begin{subfigure}{0.45\textwidth}
        \centering
        \includegraphics[width=\linewidth]{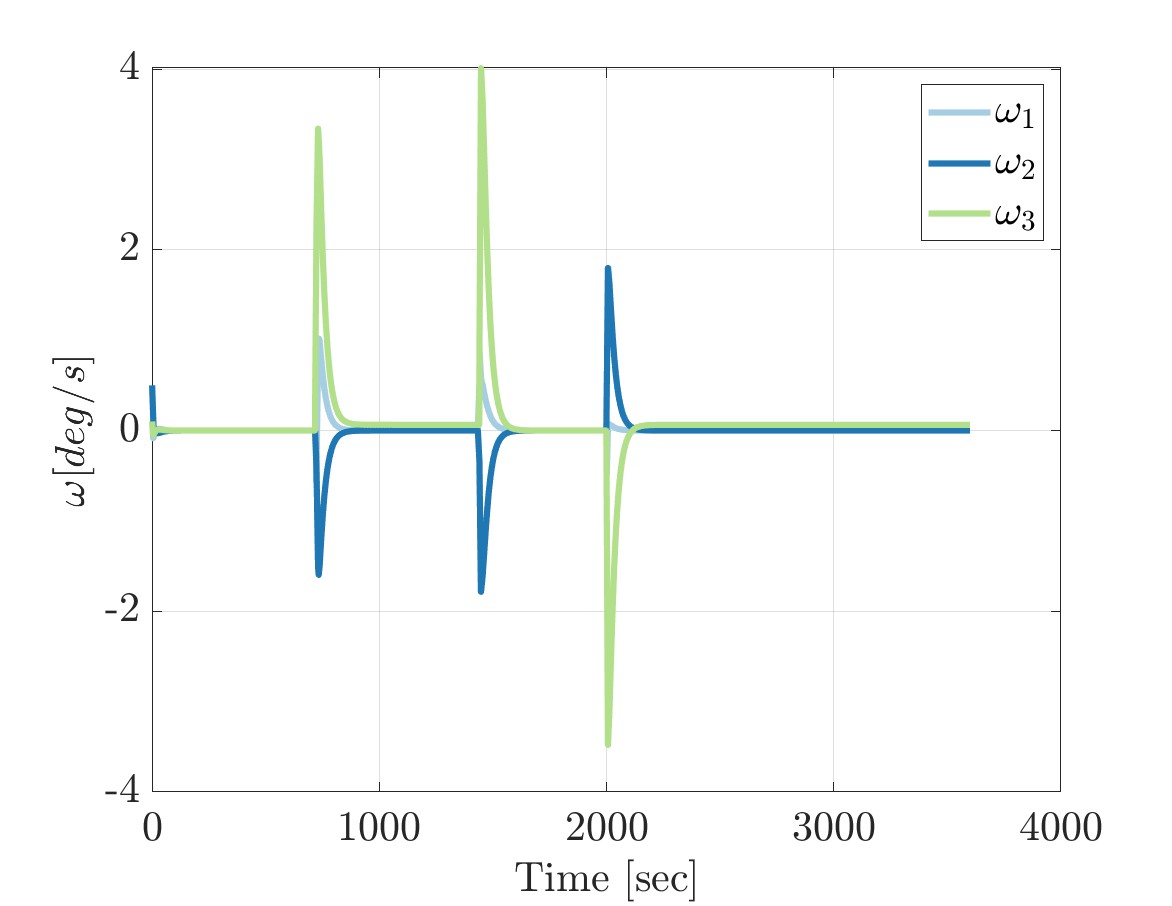}  
        \caption{Body Angular Velocity}
        \label{fig:RW4_ICL_BAV}
    \end{subfigure}

    \begin{subfigure}{\textwidth}
        \centering
        \includegraphics[width=0.45\linewidth]{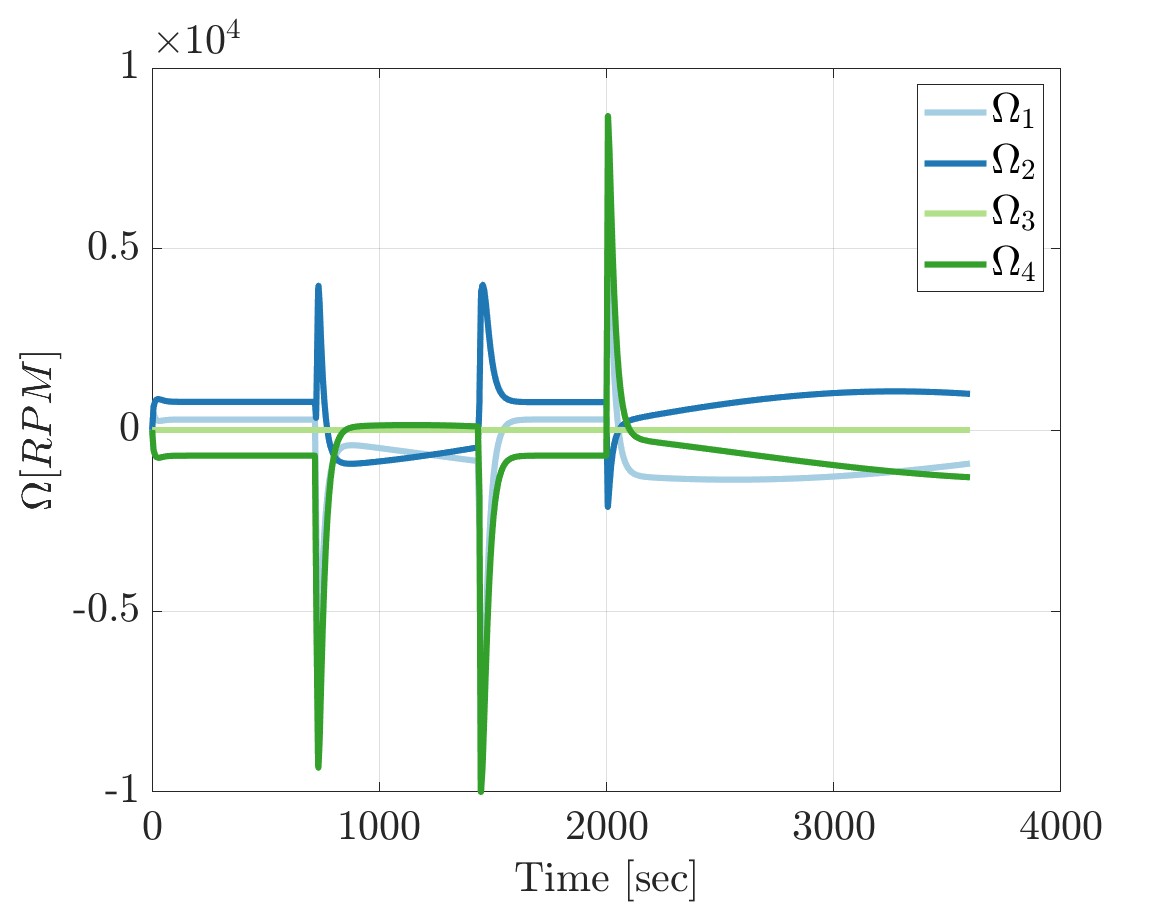}   
        \caption{Reaction Wheel Angular Velocity}
        \label{fig:RW4_ICL_AV}
    \end{subfigure}
    \caption{Case 1: 4 RW's with ICL term activated.}
\end{figure}

\begin{figure}[h]
    \begin{subfigure}{0.45\textwidth} 
        \centering
        \includegraphics[width=\linewidth]{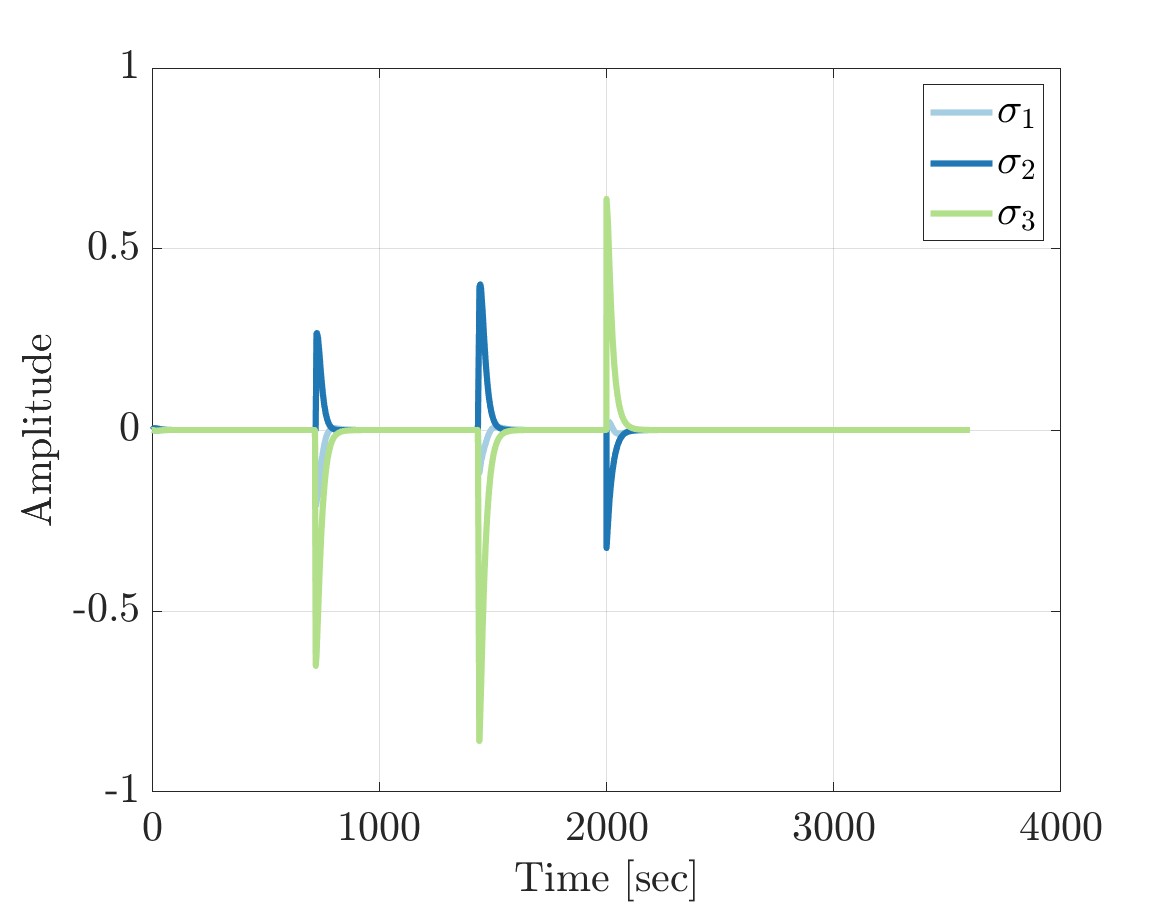}  
        \caption{Error MRP}
        \label{fig:RW4_NO_ICL_MRP}
    \end{subfigure}
    \begin{subfigure}{0.45\textwidth}
        \centering
        \includegraphics[width=\linewidth]{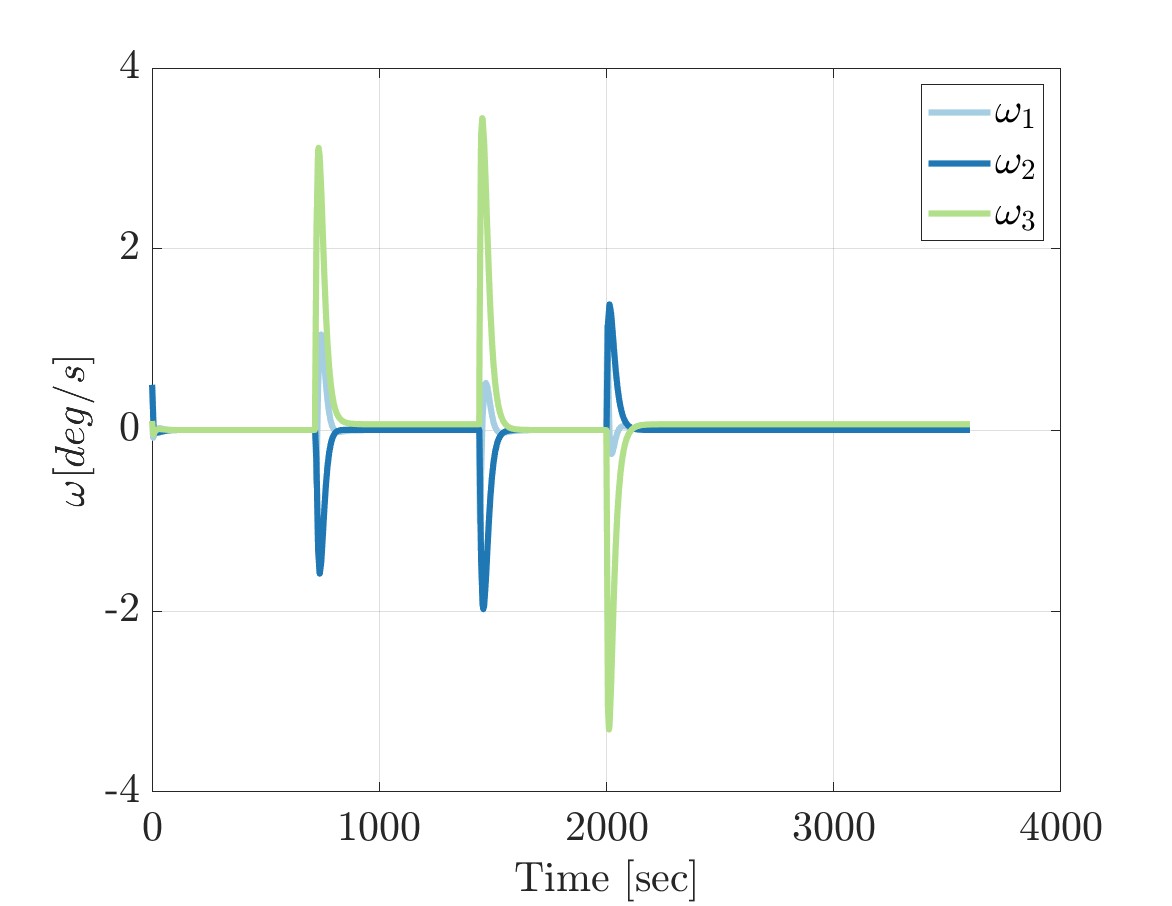}  
        \caption{Body Angular Velocity}
        \label{fig:RW4_NO_ICL_BAV}
    \end{subfigure}

    \begin{subfigure}{\textwidth}
        \centering
        \includegraphics[width=0.45\linewidth]{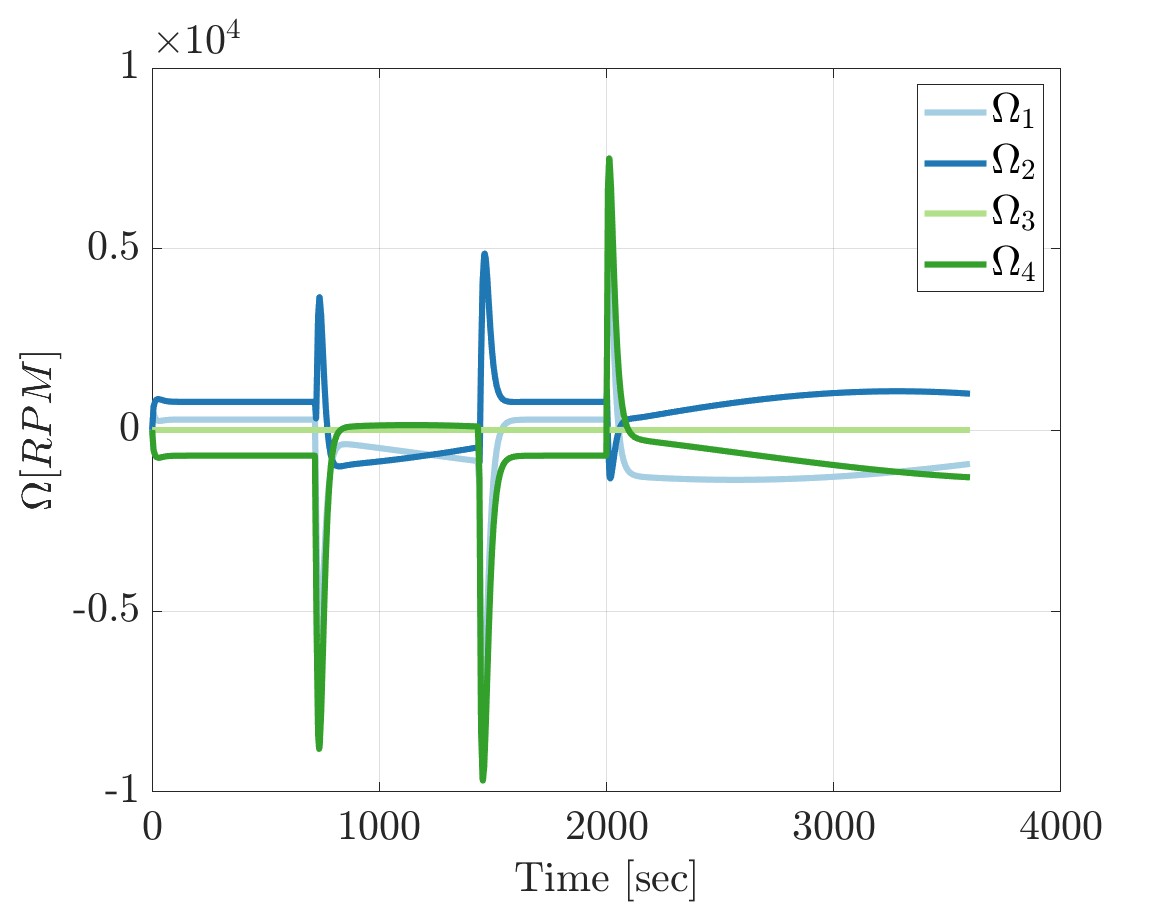}   
        \caption{Reaction Wheel Angular Velocity}
        \label{fig:RW4_NO_ICL_AV}
    \end{subfigure}
    \caption{Case 2: 4 RW's with no ICL term activated.}
\end{figure}

For both simulations, the time history plots for the error MRP, body angular velocity and RW angular velocity, are presented in  Figures \ref{fig:RW4_ICL_MRP} \&  \ref{fig:RW4_NO_ICL_MRP} ,  \ref{fig:RW4_ICL_BAV} \& \ref{fig:RW4_NO_ICL_BAV} and \ref{fig:RW4_ICL_AV} \& \ref{fig:RW4_NO_ICL_AV} for Cases 1 \& 2, respectively, to verify that the controller is accurately tracking the guidance command. The failing RW can be observed in the plot of the RWs angular velocities presented in Figures \ref{fig:RW4_ICL_AV} \& \ref{fig:RW4_NO_ICL_AV}. 

\begin{figure}[h]
    \begin{subfigure}{0.45\textwidth} 
        \centering
        \includegraphics[width=\linewidth]{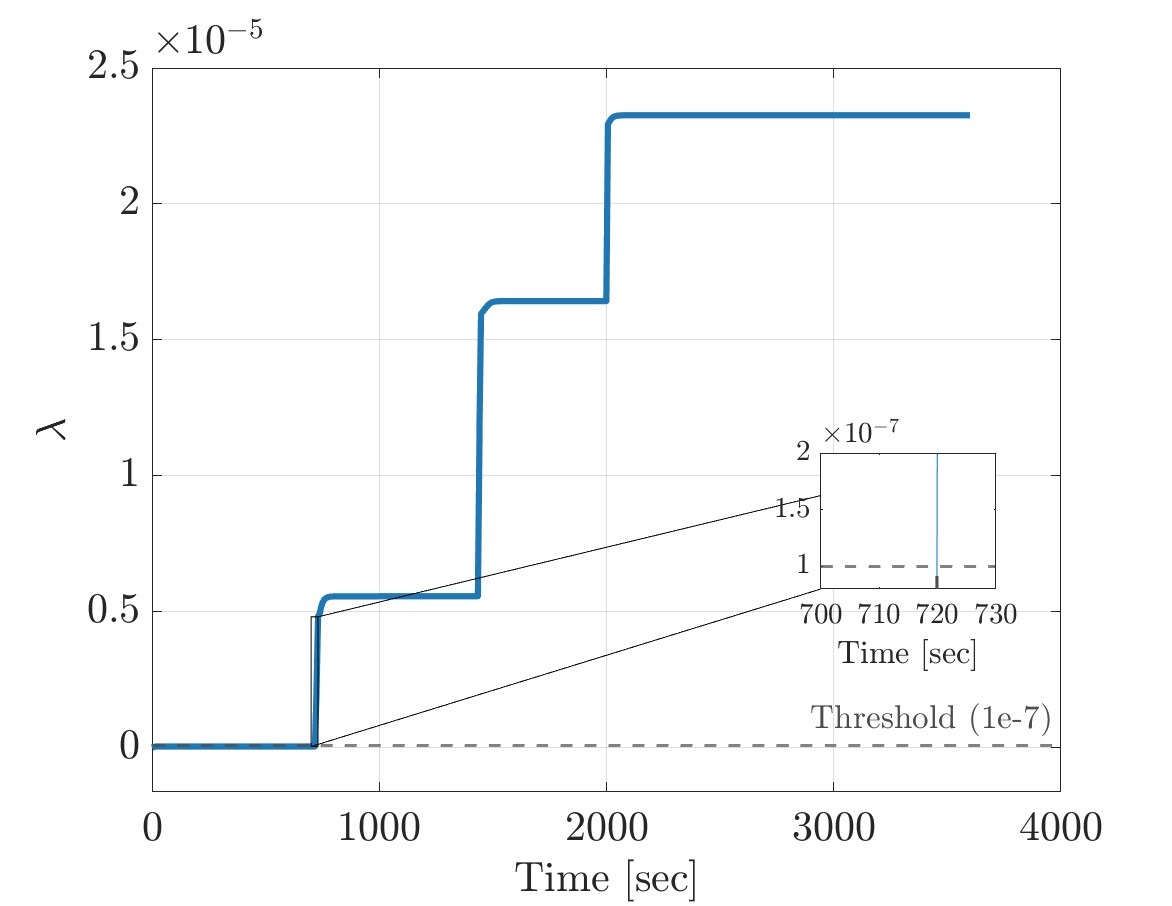}  
        \caption{Minimum Eigenvalue of Finite Excitation Condition Matrix w/ Threshold for Learning}
        \label{fig:RW4_ICL_L}
    \end{subfigure}
    \begin{subfigure}{0.45\textwidth}
        \centering
        \includegraphics[width=\linewidth]{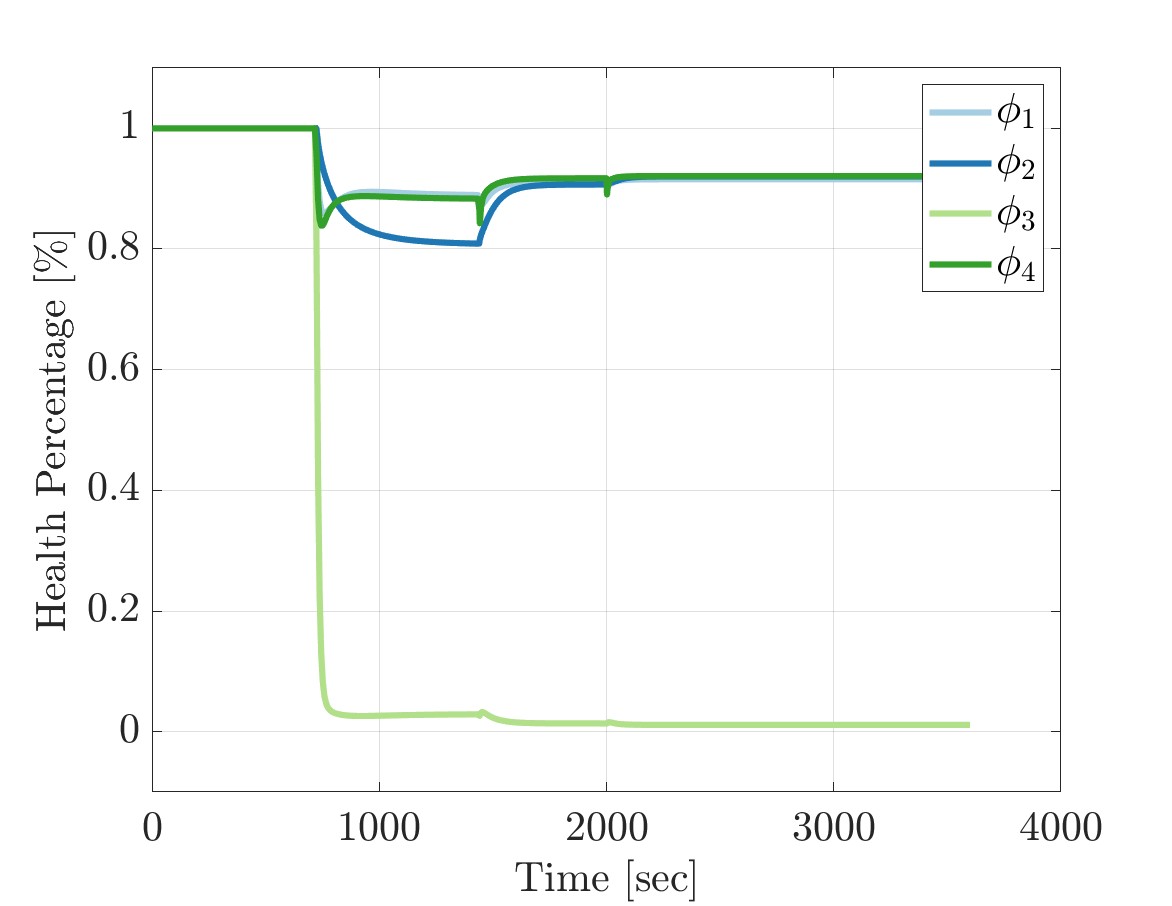}   
        \caption{Estimated Parameters}
        \label{fig:RW4_ICL_EP}
    \end{subfigure}
    \caption{Case 1: 4 RW's with ICL term activated.}
\end{figure}

\begin{figure}[h]
     \begin{subfigure}{\textwidth}
        \centering
        \includegraphics[width=0.45\linewidth]{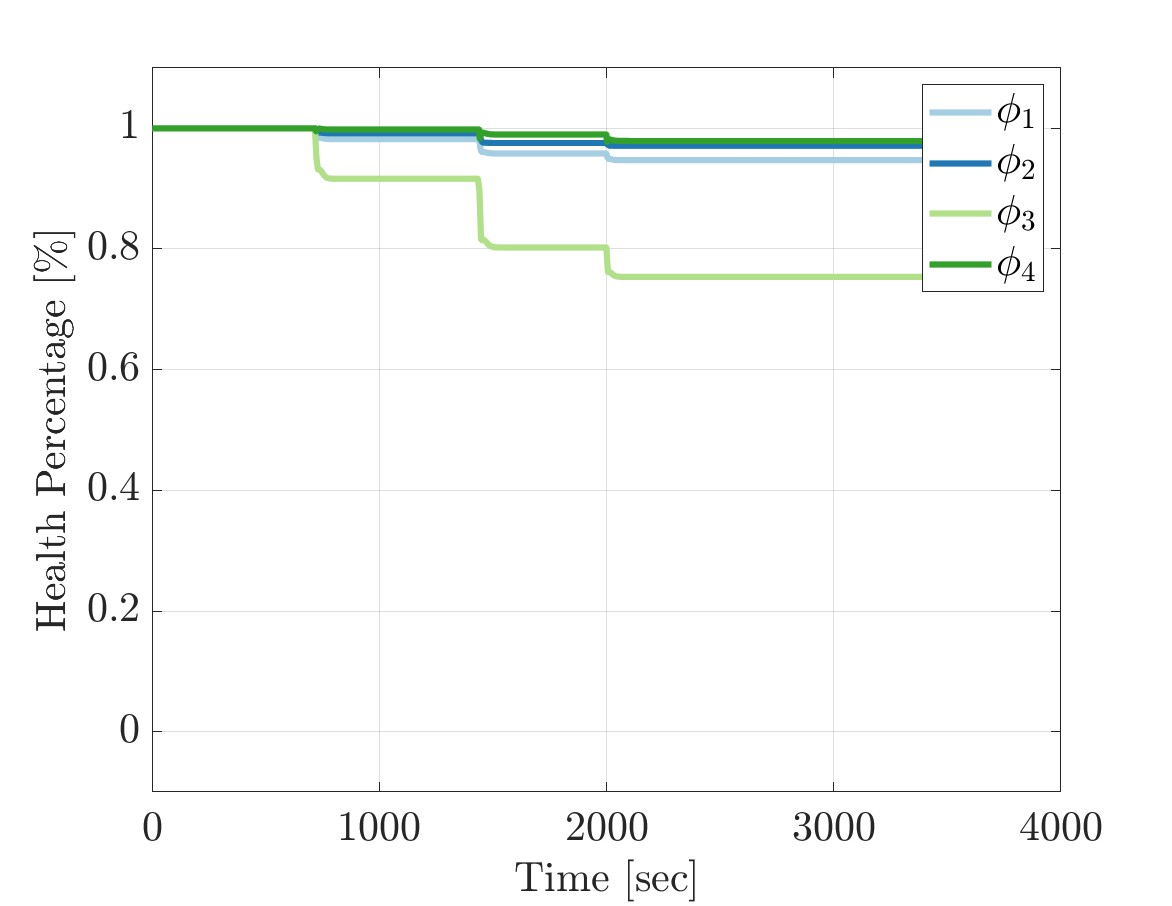}   
        \caption{Estimated Parameters}
        \label{fig:RW4_NO_ICL_EP}
    \end{subfigure}
    \caption{Case 2: 4 RW's with no ICL term activated.}
\end{figure}

For Case 1, when the ICL term is active, we include Figure \ref{fig:RW4_ICL_L}, a plot of $\lambda$, which is the minimum eigenvalue of the finite excitation condition matrix explained in Assumption 2. Because there is a user-defined threshold $\Bar{\lambda}$, which can be verified during execution, when this threshold is met it means that the system is sufficiently excited and exponential convergence of the estimation error can be guaranteed. The estimation of the uncertain parameters is presented in Figures \ref{fig:RW4_ICL_EP} \& \ref{fig:RW4_NO_ICL_EP} for both simulations. Because the ICL term is deactivated in Case 2, no plot of $\lambda$ is necessary.

\clearpage

\subsection{Cases 3 \& 4: More RWs and partial degradation}

\begin{figure}[h]
    \begin{subfigure}{0.45\textwidth} 
        \centering
        \includegraphics[width=\linewidth]{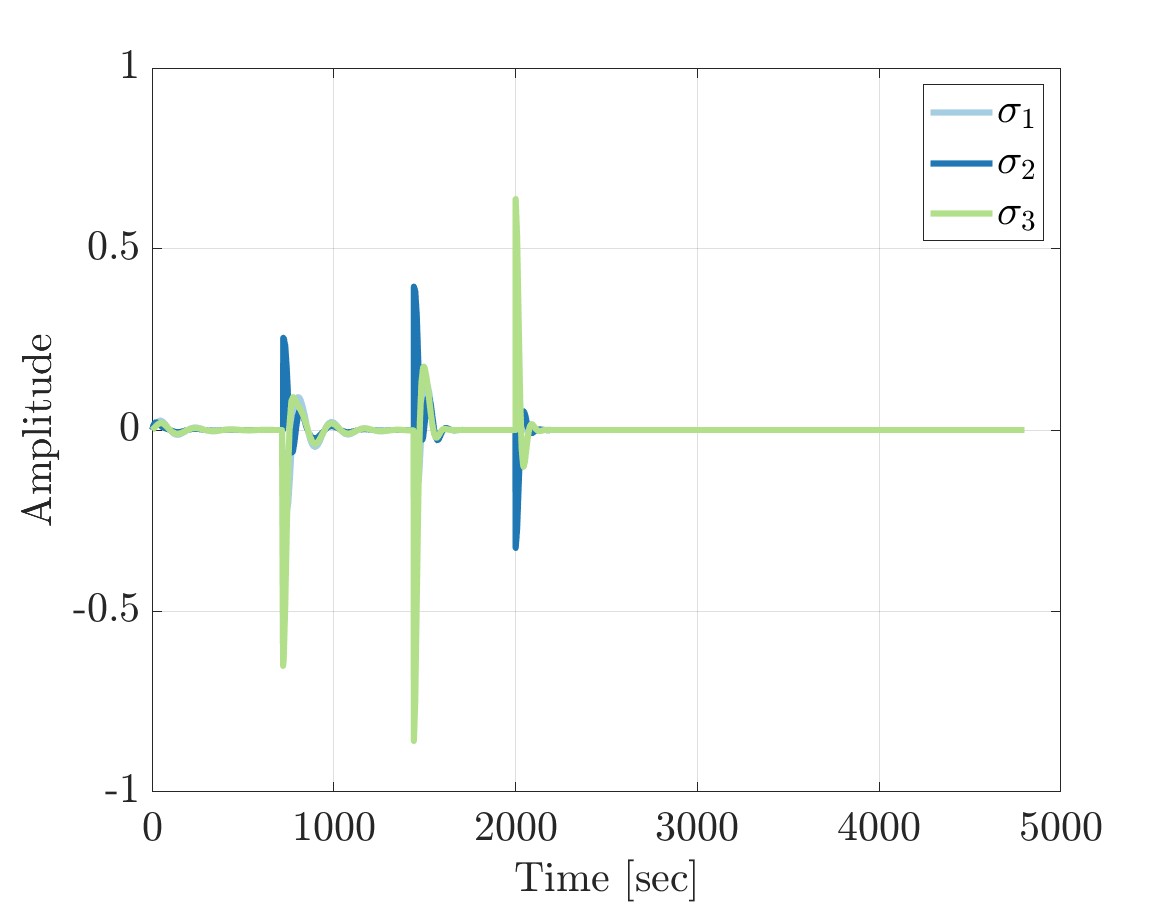}  
        \caption{Error MRP}
        \label{fig:RW6_0_MRP}
    \end{subfigure}
    \begin{subfigure}{0.45\textwidth}
        \centering
        \includegraphics[width=\linewidth]{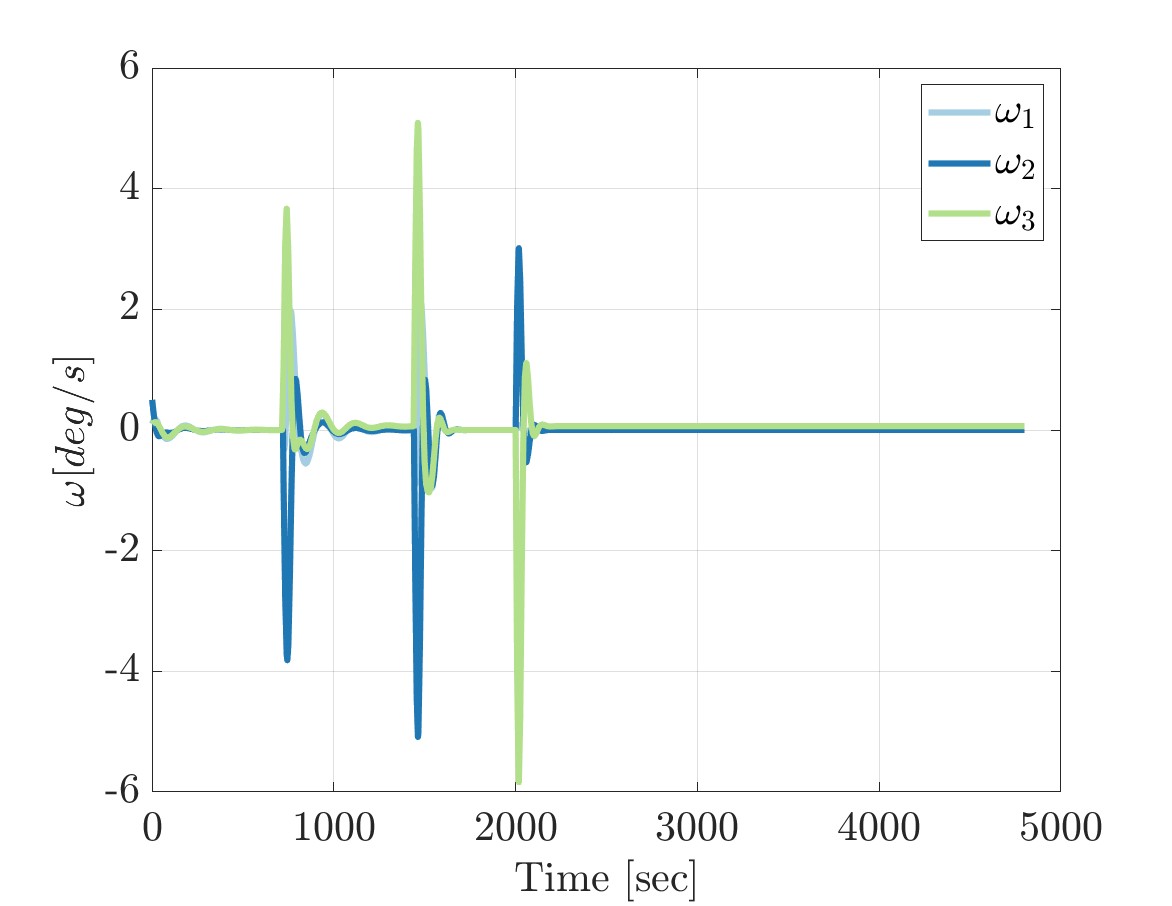}  
        \caption{Body Angular Velocity}
        \label{fig:RW6_0_BAV}
    \end{subfigure}

    \begin{subfigure}{\textwidth}
        \centering
        \includegraphics[width=0.45\linewidth]{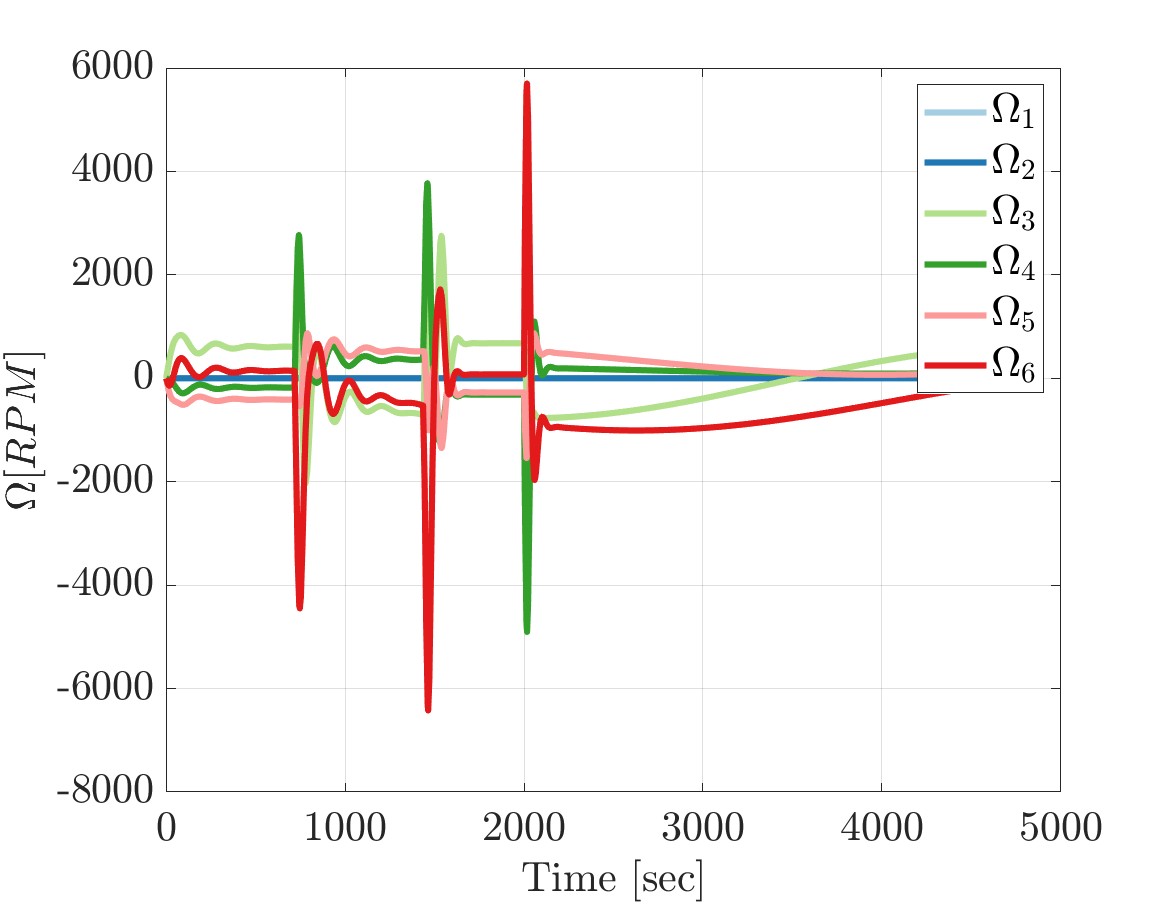}   
        \caption{Reaction Wheel Angular Velocity}
        \label{fig:RW6_0_AV}
    \end{subfigure}
    \caption{Case 3: 6 RW's with failure of RW's 1 \& 2 set to 0\%.}
\end{figure}

\begin{figure}[hp]
    \begin{subfigure}{0.45\textwidth} 
        \centering
        \includegraphics[width=\linewidth]{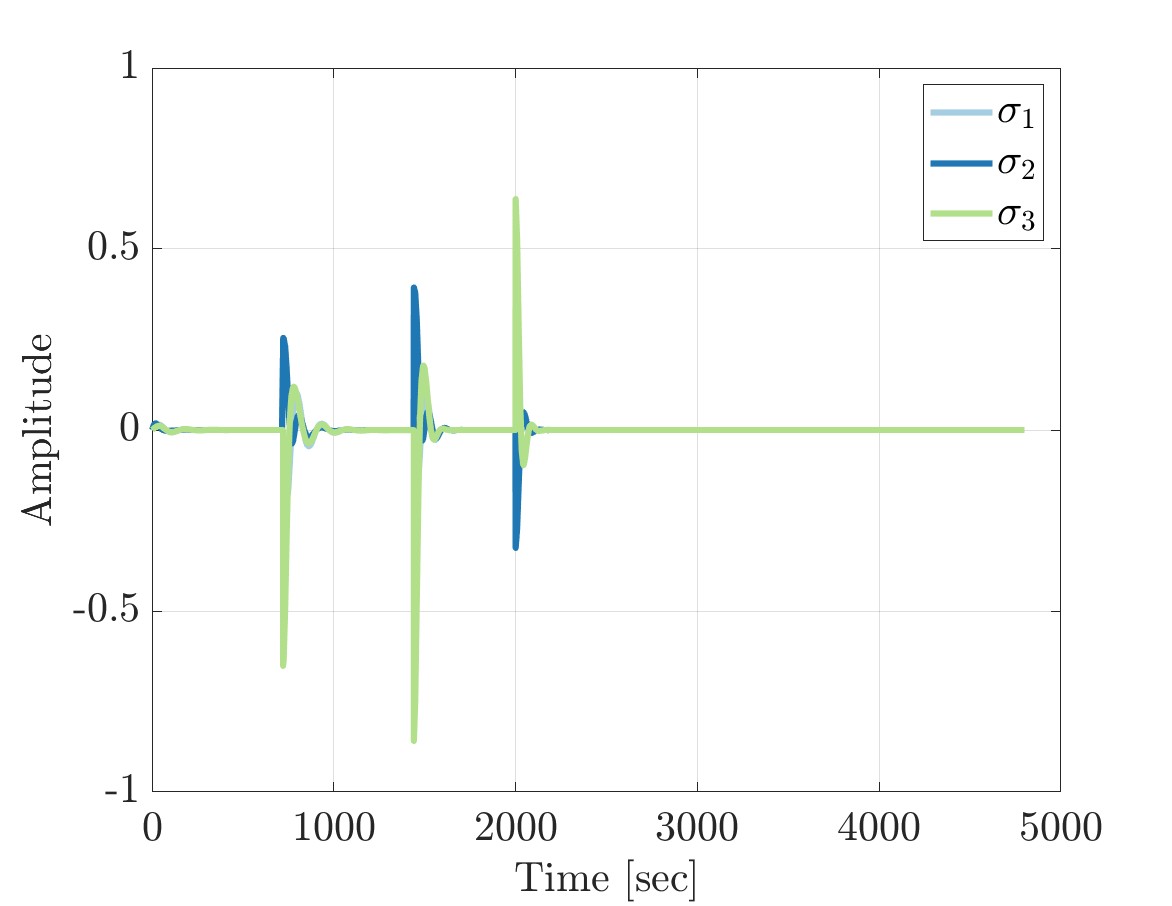}  
        \caption{Error MRP}
        \label{fig:RW6_02_MRP}
    \end{subfigure}
    \begin{subfigure}{0.45\textwidth}
        \centering
        \includegraphics[width=\linewidth]{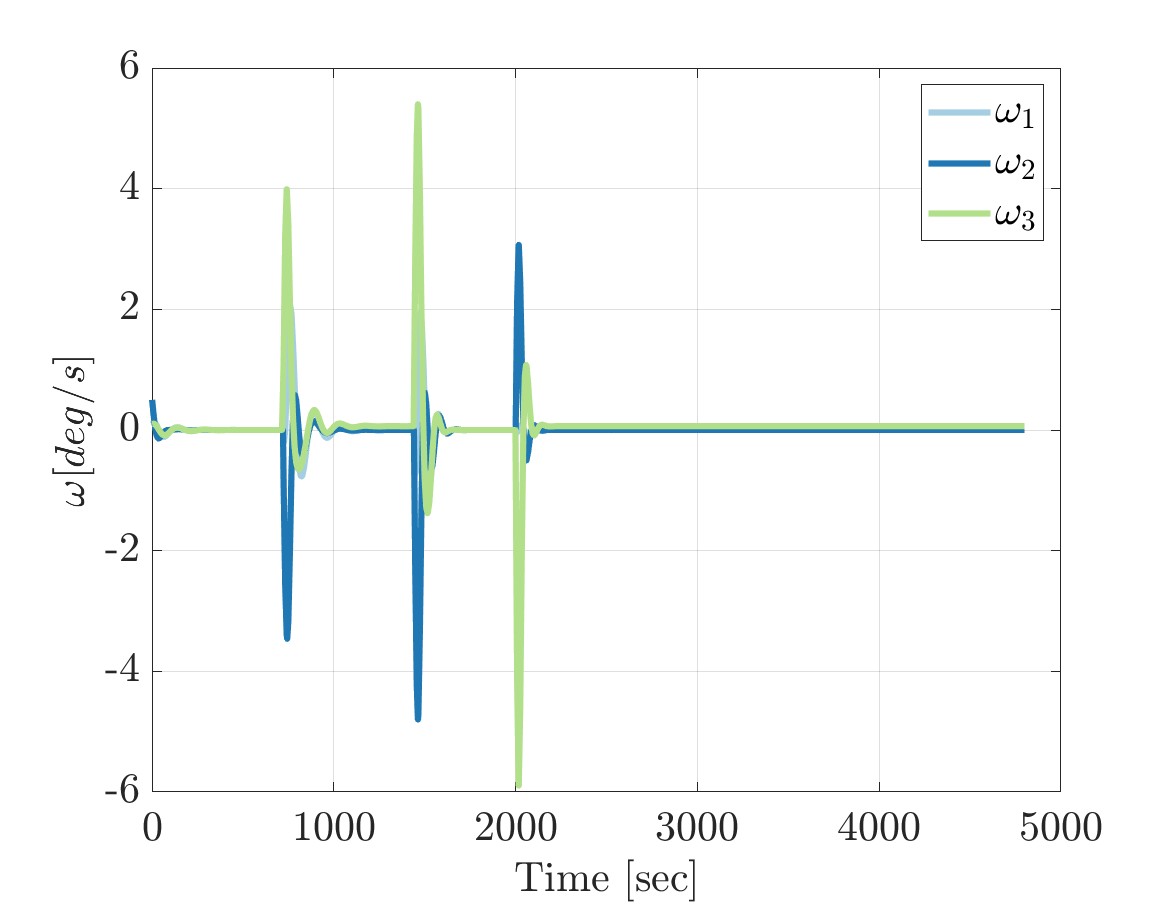}  
        \caption{Body Angular Velocity}
        \label{fig:RW6_02_BAV}
    \end{subfigure}

    \begin{subfigure}{\textwidth}
        \centering
        \includegraphics[width=0.45\linewidth]{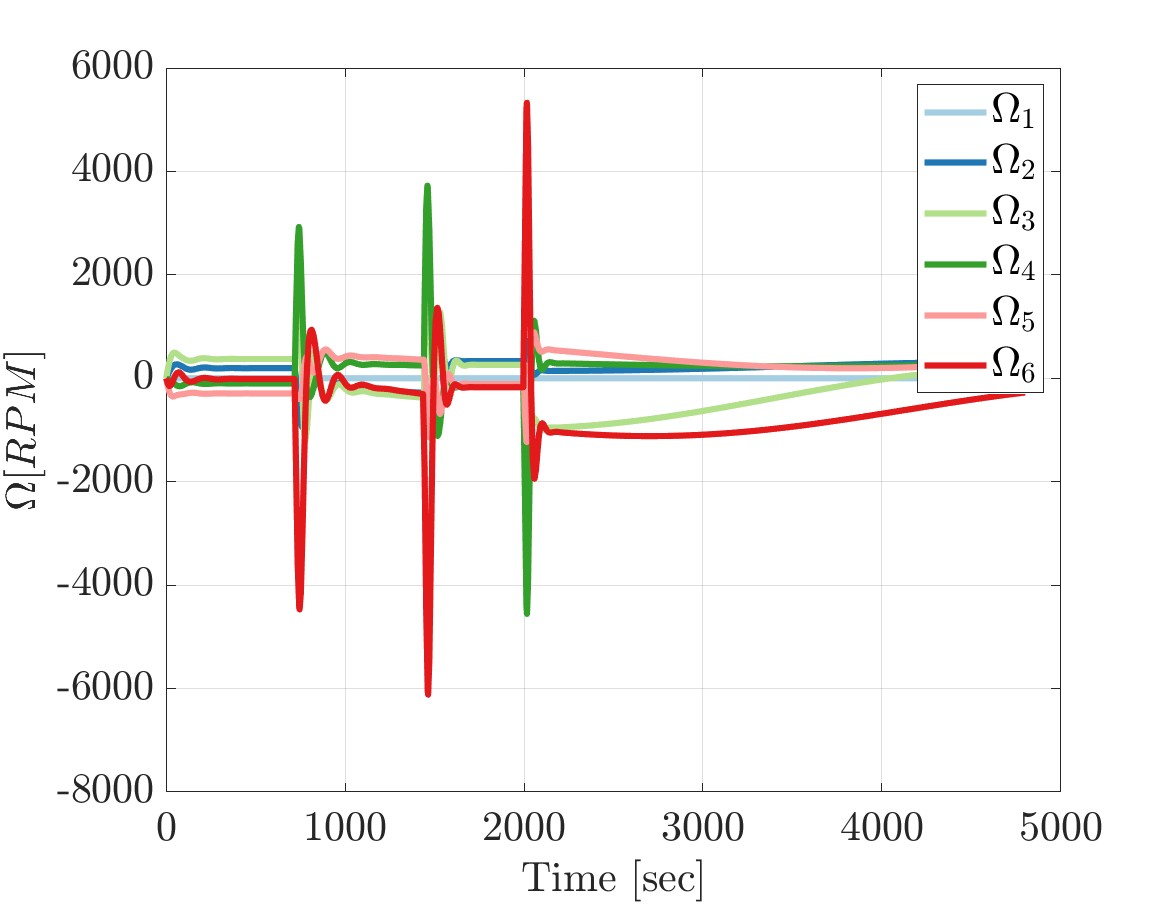}   
        \caption{Reaction Wheel Angular Velocity}
        \label{fig:RW6_02_AV}
    \end{subfigure}
    \caption{Case 4: 6 RW's with performance of RW 1 set to 0\% and RW 2 set to 30\%.}
\end{figure}

As mentioned, Cases 3 \& 4 are designed to test the capability of the proposed adaptive controller to perform in multiple scenarios and still deliver the expected results. Since Cases 1 \& 2 were designed around a system with 4 RWs, these next two simulations contain a system with 6 RWs. The 6 RWs are arranged in a vertical, symmetric, hexagon such that there are multiple redundancies built into the system. Along with the increased number of RWs, the system in both simulations will experience the failure of two RWs, to test the ability of the controller to handle multiple failures. In Case 3, RWs 1 and 2 are set to experience total failure, whilst in Case 4 the same RWs are set to experience degraded performance, only able to output 0 and 30\% of their nominal torque preformance, respectively. Like Cases 1 \& 2, for both scenarios we display the time history of their error MRP in Figures \ref{fig:RW6_0_MRP} \& \ref{fig:RW6_02_MRP}, the body angular velocity in Figures \ref{fig:RW6_0_BAV} \& \ref{fig:RW6_02_BAV} and the RW angular velocities in Figures \ref{fig:RW6_0_AV} \& \ref{fig:RW6_02_AV}. These present the overall performance of the controller and ensure that it behaves as expected. 

\begin{figure}[hp]
    \begin{subfigure}{0.45\textwidth} 
        \centering
        \includegraphics[width=\linewidth]{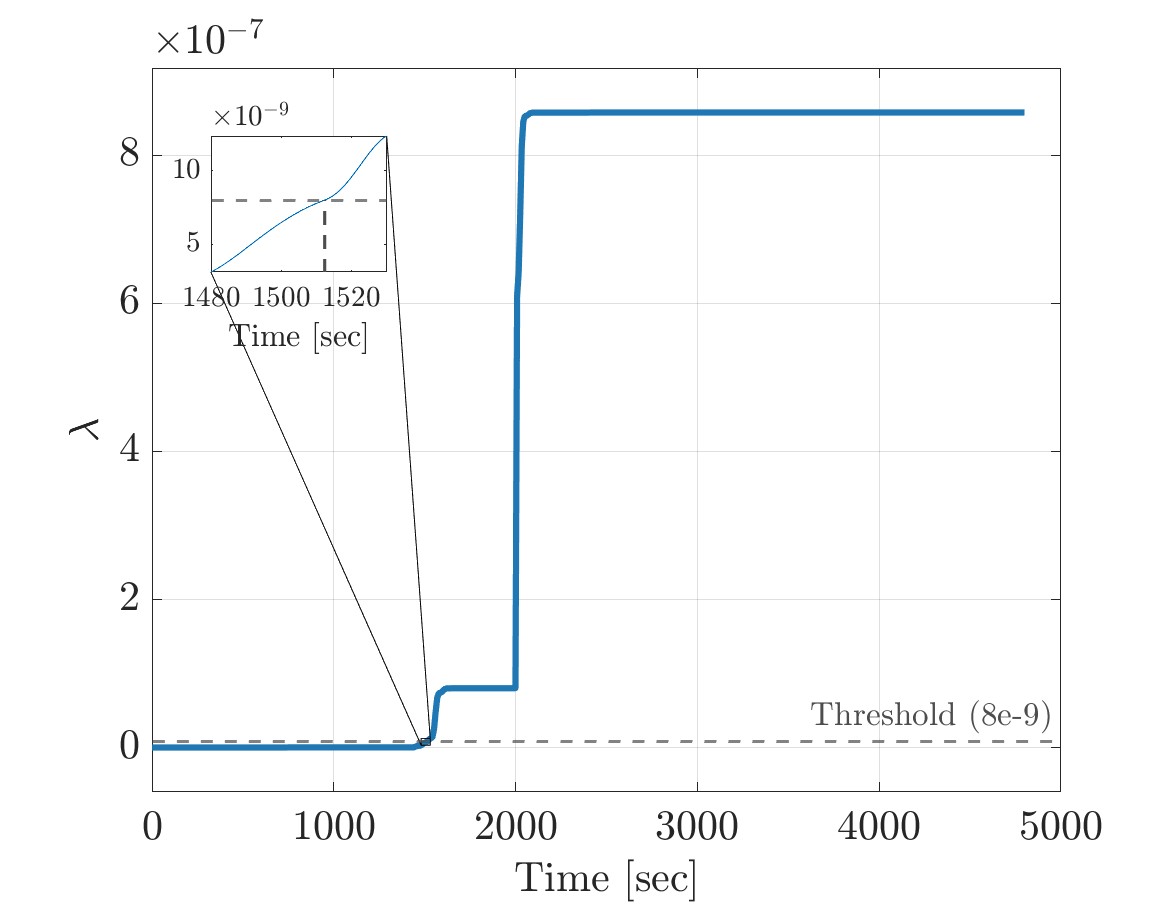}  
        \caption{Minimum Eigenvalue of Finite Excitation Condition Matrix w/ Threshold for Learning}
        \label{fig:RW6_0_L}
    \end{subfigure}
    \begin{subfigure}{0.45\textwidth}
        \centering
        \includegraphics[width=\linewidth]{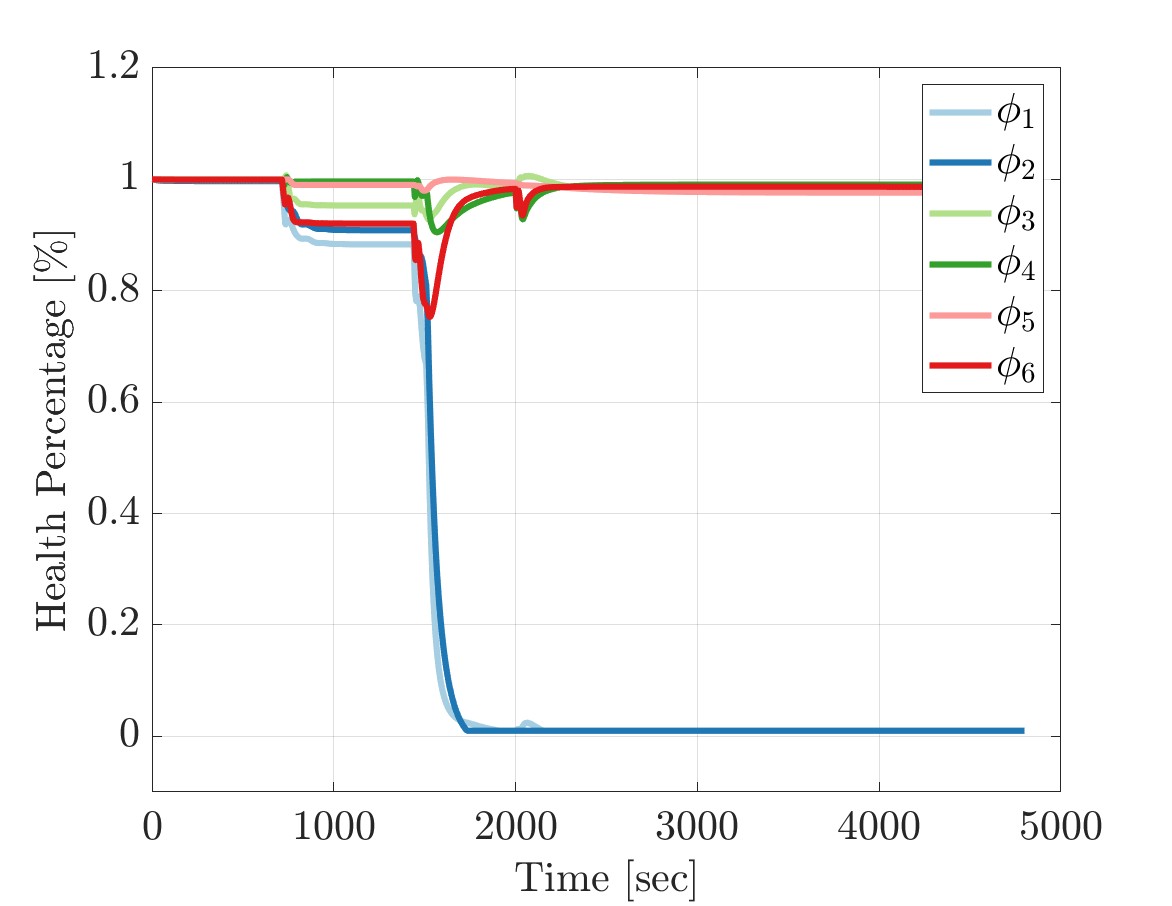}   
        \caption{Estimated Parameters}
        \label{fig:RW6_0_EP}
    \end{subfigure}
    \caption{Case 3: 6 RW's with failure of RW's 1 \& 2 set to 0\%.}
\end{figure}

\begin{figure}[htp]
    \begin{subfigure}{0.45\textwidth} 
        \centering
        \includegraphics[width=\linewidth]{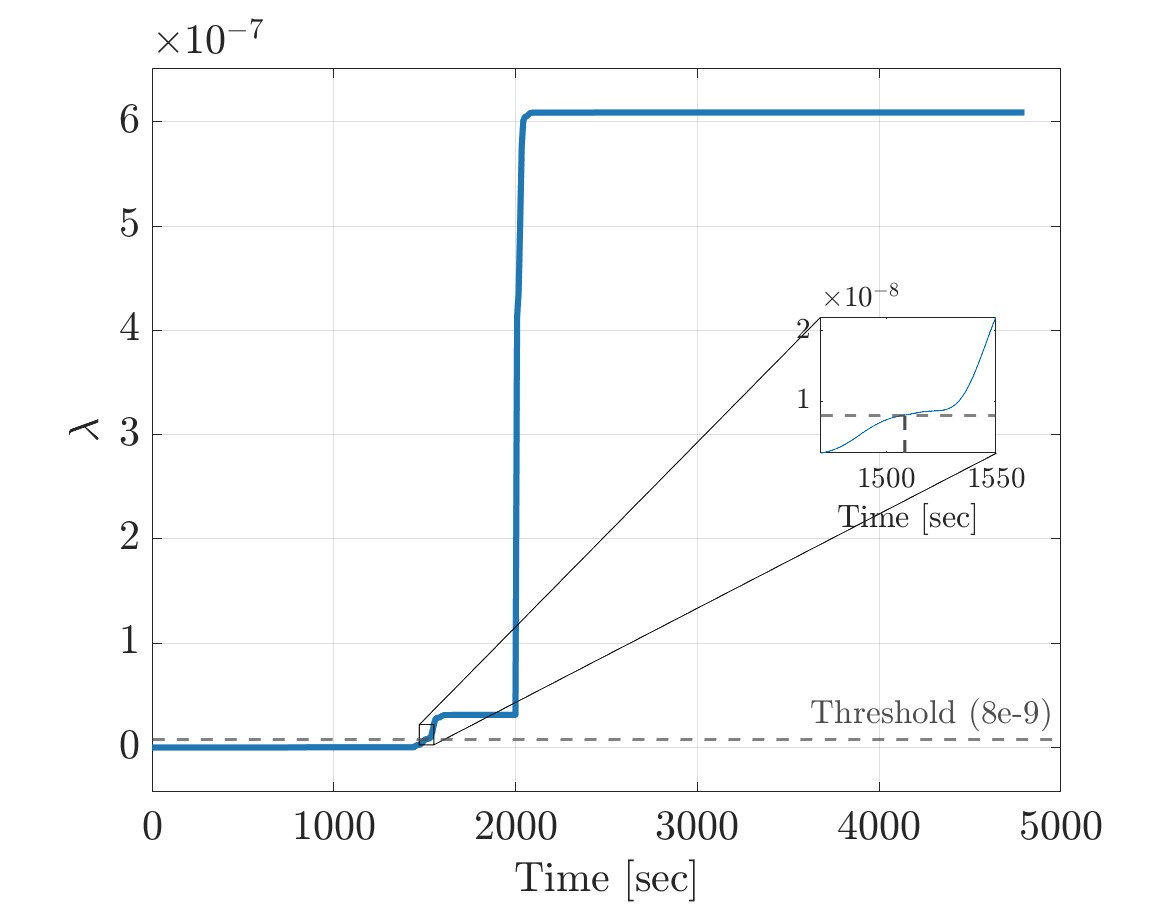}  
        \caption{Minimum Eigenvalue of Finite Excitation Condition Matrix w/ Threshold for Learning}
        \label{fig:RW6_02_L}
    \end{subfigure}
    \begin{subfigure}{0.45\textwidth}
        \centering
        \includegraphics[width=\linewidth]{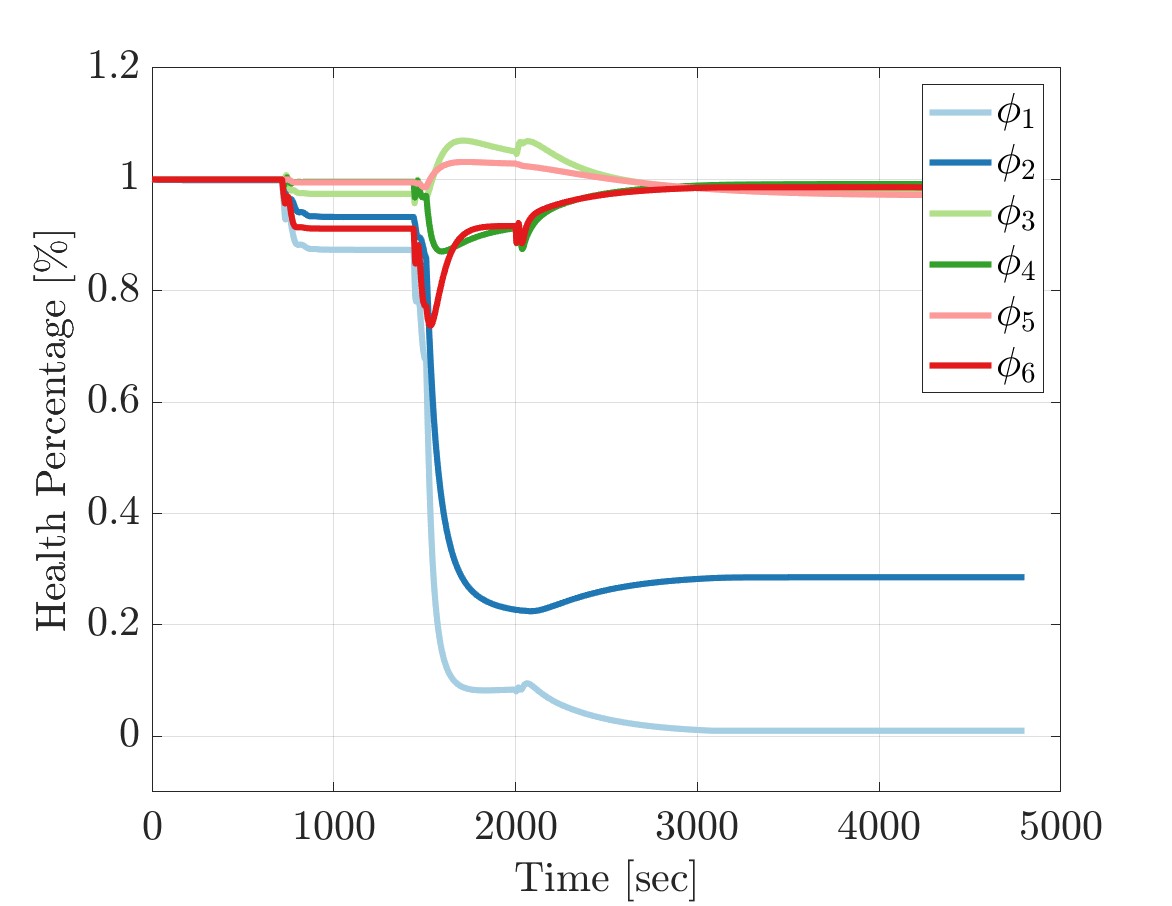}   
        \caption{Estimated Parameters}
        \label{fig:RW6_02_EP}
    \end{subfigure}
    \caption{Case 4: 6 RW's with performance of RW 1 set to 0\% and RW 2 set to 30\%.}
\end{figure}

As both simulations have the ICL term active, Figures \ref{fig:RW6_0_L} \& \ref{fig:RW6_02_L} display the value of $\lambda$ over time, whilst Figures \ref{fig:RW6_0_EP} \& \ref{fig:RW6_02_EP} show the values of the estimated uncertain RW health parameters. 

\subsection{Results}
The simulations performed in Cases 1 \& 2 illustrate the importance of the ICL term. Looking at Figures \ref{fig:RW4_ICL_EP} and \ref{fig:RW4_NO_ICL_EP} it can be observed that when the ICL term is activated, the controller is able to effectively estimate the correct health value for all RWs, specifically the fact that RW 3 loses all of its effective performance. In the simulataion with no ICL term, the estimated RW health parameters in $\boldsymbol{\hat{\theta}}$ do not converge to their real values. This is expected because the stability result from theorem 1 does not ensure convergence of the estimation error $\boldsymbol{\tilde{\theta}}$. The periods of increased growth in $\lambda$ at  720, 1440, 2000 seconds are indicative of the system being excited due to the change in pointing reference. During these periods where the RWs are being more active to alter the pointing of the satellite, the system is able to gather input-output data through the ICL term until the threshold is reached and accurate estimation is guaranteed thereafter. This is evident in Figure \ref{fig:RW4_ICL_EP} as the exponential decay to 0\% health occurs when the threshold is reached, as expected from theorem 2. 

In both scenarios, it is important to note that the controller was still able to correctly control the satellite to achieve accurate tracking of the attitude guidance. However, without the ICL term, the controller was unable to learn the uncertain parameters of the system. In both scenarios, the remaining RWs did not saturate with the failure of the third RW.

Cases 3 \& 4 investigate the capabilities of the controller to handle a more complex system with an increased number of RWs as well as multiple failures at different levels. In Case 3, both RWs 1 and 2 were set to completely fail and this was correctly detected by the adaptive controller as seen in Figure \ref{fig:RW6_0_EP}. Once the threshold was crossed, and the system was sufficiently excited, the estimation of the health of the RWs was correctly estimated. Likewise, in Case 4, where the RWs 1 and 2 were set to be degraded to 0 and 30\% performance, respectively, the controller was able to accurately determine the health and estimate these uncertainties. An important note for both Case 3 \& 4, is that although two of the six RWs failed, the other four RWs did not saturate, a critical factor in satellite attitude control, showing that no excessive control effort was added once the FE condition was satisfied. 

As in Cases 1 \& 2, both simulations were able to accurately track the attitude guidance, showing that despite multiple failures, the system was still able to perform its control task and provide correct estimations.

\begin{figure}[hp]
    \begin{subfigure}[b]{0.45\textwidth} 
        \centering
        \includegraphics[width=\textwidth]{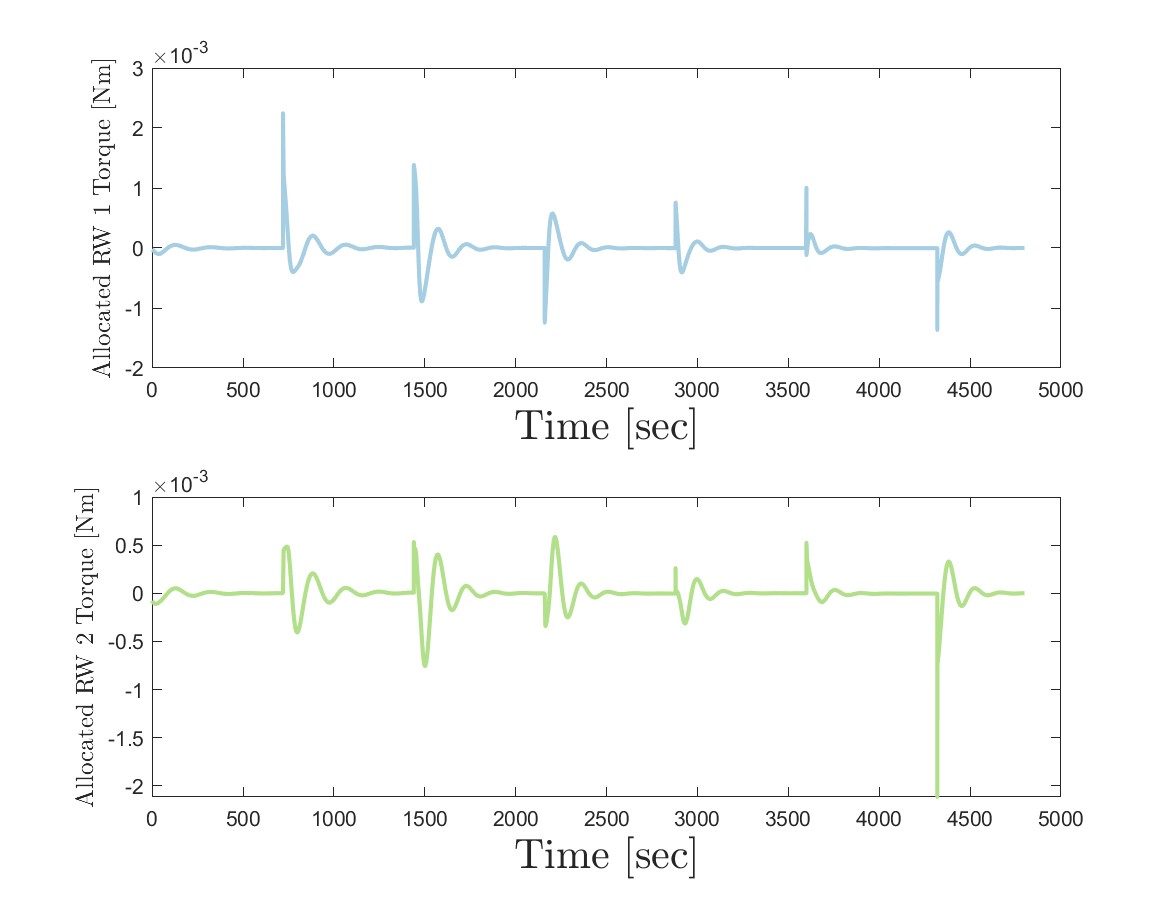}  
        \caption{Controller allocated torque to RW's 1 and 2 for the simulation with no ICL term active.}
        \label{fig:RW6_NO_ICL_0_T}
    \end{subfigure}
    \hfill
    \begin{subfigure}[b]{0.45\textwidth}
        \centering
        \includegraphics[width=\textwidth]{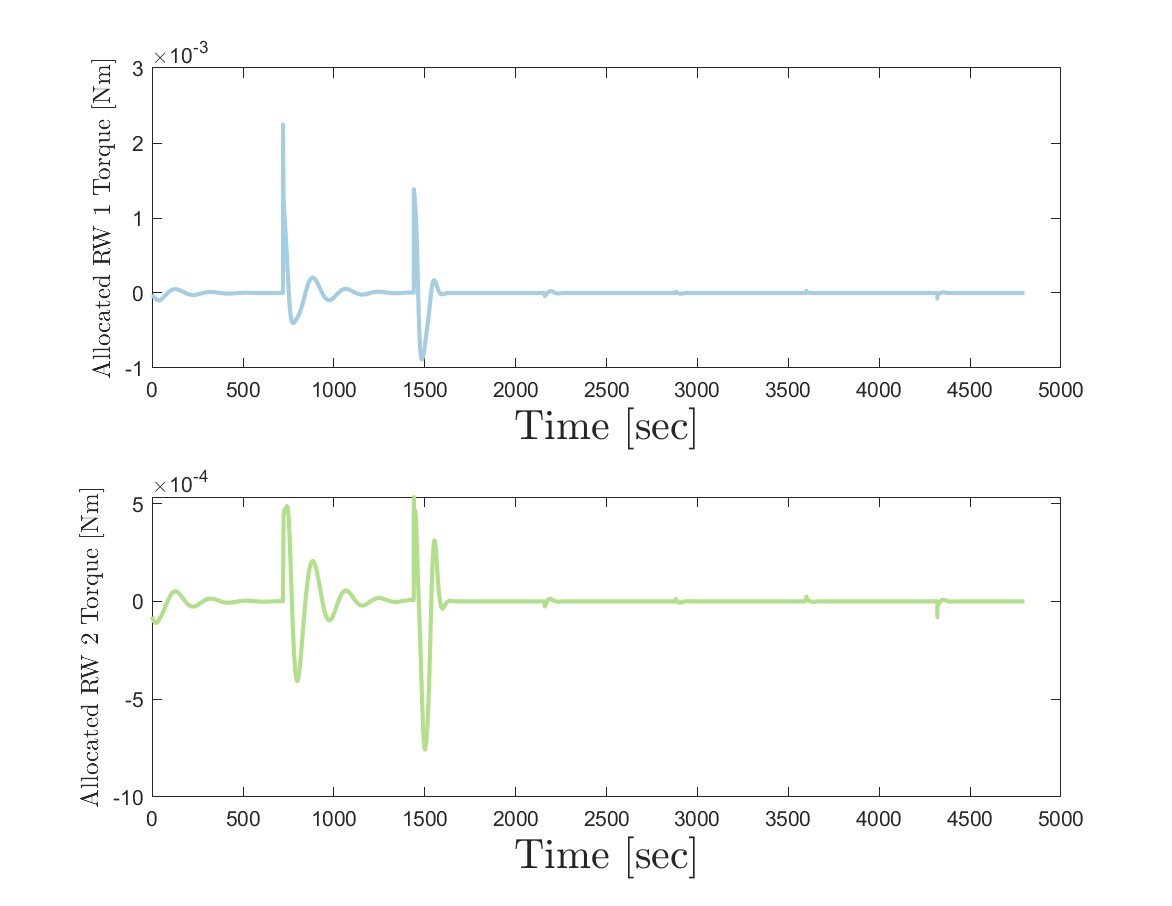}   
        \caption{Controller allocated torque to RW's 1 and 2 for the simulation with the ICL term active.}
        \label{fig:RW6_0_T}
    \end{subfigure}
    \caption{Comparison of allocated torques on the degraded RW's when the ICL term is active or not.}
\end{figure}

An interesting observation of the performance of the adaptive controller can be noted when running the same simulation, with and without the ICL term. Figures \ref{fig:RW6_NO_ICL_0_T} and \ref{fig:RW6_0_T} show the torque allocated by the adaptive controller to both RW 1 and 2 in each scenario with 6 total RW's. In both simulations, all parameters, but the ICL term, were kept equal, whilst RW 1 and 2 were degraded to 0\% performance. The key difference between these simulations is that in Figure \ref{fig:RW6_NO_ICL_0_T} the ICL term is not active, whereas in Figure \ref{fig:RW6_0_T} the ICL term is active. We can clearly see that when the ICL term is active, and the learning threshold is reached ($\sim$ 1500sec), the controller has correctly estimated that RW 1 and 2 are completely degraded, hence does not allocate any control torque to them. Figure \ref{fig:RW6_NO_ICL_0_T} shows the contrary, where the even though the RW's are degraded, the adaptive controller does not learn the estimated health parameters and continues trying to allocate control torques to these reaction wheels. 

\subsection{Discussion}
The performance of the adaptive controller in all scenarios was the expected. The ICL term was shown to be instrumental to learn the uncertain RW health parameters, whilst the controller was shown to adapt to compensate a variety of RW failures to different levels. One of the key issues in the simulation scenarios occurred in Case 1, where the estimation of the health of the non-degraded RWs is not as accurate. This could be explained by a number of reasons, primarily because the simulations run include external disturbances and nonlinearities not accounted for in the design of the controller. A refinement of the design along with further tuning of the key parameters outlined in Table \ref{tab:gains} could also improve the performance of the estimation. We see however in Cases 3 \& 4, excellent performance of the estimation for all RWs, with the working reaction wheels all being estimated very close to 100\% health whereas the degraded RWs reach a steady state value very close to their true health. 

\section{Conclusion}

In this work, an ICL-based adaptive controller that can simultaneously estimate and compensate for the degradation of reaction wheels was designed. Through the implementation of an adaptive update law that incorporates integral concurrent learning, the controller successfully estimated the uncertain parameters, namely the RWs health. The finite excitation condition was verified online and the uncertain parameters estimation converged exponentially afterwards, as expected from the presented Lyapunov-based stability analysis. Posterior to the excitation condition being reached, the adaptive controller shows that it uses the learned estimated health parameters to better allocate control torques on those wheels still operable and not to those with degraded performance.

Four different scenarios were performed and the controller was able to ensure attitude tracking without saturation of the RWs, as well as estimation of the degradation level for all RWs. Although the controller was able to operate to auspiciously achieve global exponential tracking despite the multiple RW failures induced in simulation, future work could include the integration of robustness to additive disturbances in this type of controller, since these are expected to affect the stability result for both state and estimation errors. 

\bibliographystyle{AAS_publication}   
\bibliography{references}   

\end{document}